\documentclass[10pt,a4paper,twoside]{article} 
\usepackage{amsmath}
\usepackage{amssymb} \usepackage{latexsym}
\usepackage{amscd}
\usepackage{fancyhdr}
\usepackage{mathrsfs} 
\usepackage{texdraw}

\newtheorem{definition}{\bf Definition}[section]

\newtheorem{lemma}[definition]{\bf Lemma}
\newtheorem{theorem}[definition]{\bf Theorem}
\newtheorem{proposition}[definition]{\bf Proposition}
\newtheorem{conjecture}[definition]{\bf Conjecture}

\newcounter{remark}
\newcounter{remark:int}

\newenvironment{remark}{\setcounter{remark}{\value{definition}}
  \addtocounter{remark}{1} \addtocounter{definition}{1}
  \noindent {\bf Remark \arabic{section}.\arabic{remark}} }

\newenvironment{acknowledgements}{\medskip\noindent\footnotesize{\bf
    Acknowledgements:}}{\normalsize}

\newcommand{\qed}{\rightline{$\Box $}}
\newcommand{\wee}{{\scriptscriptstyle \vee}}
\newcommand{\Db}{\mathfrak{Db}}
\newcommand{\Dbrd}[1]{\mathfrak{Db}_{\widetilde{X}}^{{\rm rd}, #1}}
\newcommand{\Crdt}{{\cal C}^{\rm rd}_{\widetilde{X}}}
\newcommand{\Crdtp}[1]{{\cal C}^{{\rm rd}, #1}_{\widetilde{X}}}

\newcommand{\widebar}[1]{\overline{#1}} 
\newcommand{\setmin}{\! \smallsetminus \!}
\newcommand{\R}{{\mathbb R}} \newcommand{\C}{{\mathbb C}}
\newcommand{\N}{{\mathbb N}} \newcommand{\Q}{{\mathbb Q}}
\newcommand{\Ct}{{\cal C}_{\widetilde{X}}} 
\newcommand{\Ctd}{{\cal C}_{\widetilde{X}, \widetilde{D}}}
\newcommand{\Z}{{\mathbb Z}} 
\newcommand{\Hy}{{\mathbb H}}
\newcommand{\Hom}{{\rm Hom}} 
\newcommand{\cHom}{\mathcal{H}om}

\newcommand{\Epq}[2]{{\Omega^{\infty, (#1,#2)}_{\Xan}}}

\newcommand{\Etop}[1]{{\Omega^{\infty,#1}_{\widetilde{X}}}}

\newcommand{\EYpq}[2]{{\Omega^{\infty, (#1,#2)}_{Y^{\rm an}}}}
\newcommand{\MM}{{\cal M}}
\newcommand{\W}{{\cal W}}

\newcommand{\vi}{\varphi} \newcommand{\ve}{\varepsilon}
\newcommand{\bve}{{\cal E}} \newcommand{\vt}{\vartheta}
 \newcommand{\Int}{\int\limits}
 
\newcommand{\lto}{\longrightarrow} 
\newcommand{\Crd}{{\cal C}^{\rm rd}} 
\newcommand{\Hrd}{H^{rd}} 
\newcommand{\Hdr}{H_{dR}} 

\newcommand{\M}{{\cal M}}
 \renewcommand{\H}{\mathcal{H}}
\renewcommand{\O}{{\cal O}} \newcommand{\pf}{{\bf Proof: }}
\newcommand{\OXY}{\O_{\widehat{X|Y}}} 
\newcommand{\A}{{\cal A}_{\widetilde{X}}} 
\newcommand{\AD}{{\cal A}_{\widetilde{X}}^{<D}} 
\newcommand{\fA}{{\cal A}_{\widehat{\widetilde{X} | D}}}
\newcommand{\dr}[2]{{\rm DR}_{\widetilde{#1}}^{#2}}
\newcommand{\DR}{{\rm DR}_{\Xan}} 
 
\newcommand{\DRD}{{\rm DR}^{<D}_{\widetilde{X}}}

\newcommand{\Amod}{{\cal A}_{\widetilde{X}}^{{\rm mod }D}} 
\newcommand{\DRmod}{{\rm DR}_{\widetilde{X}}^{{\rm mod }D}}
\newcommand{\DRmd}[1]{{\rm DR}_{\widetilde{X}}^{{\rm mod }D,#1}}
\newcommand{\DRd}[1]{{\rm DR}_{\widetilde{X}}^{<D,#1}}
\newcommand{\DRmY}{{\rm DR}_{\widetilde{Y}}^{{\rm mod }S}}
\newcommand{\DRdY}{{\rm DR}_{\widetilde{Y}}^{<S}}
\newcommand{\DRmvt}[1]{{\rm DR}_{\widetilde{X},#1}^{{\rm mod }D}}

\newcommand{\bR}{\boldsymbol{R}}
\newcommand{\As}[2]{{\cal A}_{\widetilde{#1}}^{#2}}
\newcommand{\Avt}[1]{{\cal A}_{\widetilde{X},#1}}

\newcommand{\Amodvt}[1]{{\cal A}_{\widetilde{X},#1}^{{\rm mod }D}}
\newcommand{\AmdY}{{\cal A}_{\widetilde{Y}}^{{\rm mod }S}}

\newcommand{\U}{{\cal U}}
\newcommand{\Fmod}{{\cal F}^{{\rm mod}}}

\renewcommand{\Re}{{\rm Re} \,}
\newcommand{\Lotimes}{\overset{\mathbb L}{\otimes}}
\newcommand{\Smod}{\mathcal{S}^{{\rm mod }D}}
\newcommand{\SD}{{}^\wee\!{\mathcal{S}}^{<D}}
\newcommand{\SDo}{{\mathcal{S}}^{<D}}
\newcommand{\dolb}{\widebar{\partial}}
\newcommand{\Emd}{{\cal P}^{{\rm mod} D}}

\newcommand{\bdot}{\raisebox{-.3ex}{\mbox{\boldmath $\cdot $}}}

\newcommand{\CrdY}[1]{{\cal C}_{#1}^{{\rm rd}}}
\newcommand{\Rd}{{\cal R}d}

\newcommand{\Uan}{U^{\rm an}}
\newcommand{\Xan}{X^{\rm an}}
\newcommand{\an}{{\rm an}}

\makeatletter \renewcommand{\section}{\@startsection{section}{1}{\z@}%
  {-\baselineskip}{1\baselineskip}{\large\bf}}
\renewcommand{\subsection}{\@startsection{subsection}{2}{\z@}%
  {-\baselineskip}{0.5\baselineskip}{\normalsize\bf}}
\renewcommand{\subsubsection}{\@startsection{subsubsection}{3}{\z@}%
  {-\baselineskip}{0.5\baselineskip}{\normalsize\bf}} \makeatother

\numberwithin{equation}{section}

\begin{document}

\pagestyle{fancy} \thispagestyle{plain} \headheight12pt
\renewcommand{\headrulewidth}{0.1pt} \fancyhf{} \fancyhead[ER]{Periods
  for irregular singular connections on surfaces}
\fancyhead[EL]{\thepage} \fancyhead[OL]{Periods for irregular singular
  connections on surfaces} \fancyhead[OR]{\thepage}

\centerline{\Large\bf Periods for irregular singular connections}
\centerline{\Large\bf on surfaces}

\bigskip

\centerline{{\bf Marco Hien}}

\bigskip\noindent 
\centerline{\parbox[b]{11cm}{{\bf Abstract:} Given an integrable
    connection on a smooth quasi-projective algebraic surface $U $
    over a subfield $k $ of the complex
    numbers, we define {\em rapid decay homology} groups with respect
    to the associated analytic connection which pair with the
    algebraic de Rham cohomology in terms of period integrals. These
    homology groups generalize the analogous groups in the same
    situation over curves defined by S.~Bloch and H.~Esnault. In
    dimension two, however, new features appear in this context which
    we explain in detail. Assuming a conjecture of C.~Sabbah on the
    formal classification of meromorphic connections on surfaces
    (known to be true if the rank is lower than or equal to 5), we prove
    perfectness of the period pairing in dimension two.}} 

\section{Introduction}

\pagestyle{fancy}

Given a smooth $n $-dimensional algebraic variety $U $ over a subfield
$k \subset \C $, the algebraic de Rham cohomology $\Hdr^p(U) $ is
defined to be the hypercohomology of the de Rham complex  $0 \to \O_U
\to \Omega^1_{U|k} \to \ldots \to \Omega^n_{U|k} \to 0 $
consisting of the sheaves of K\"ahler differentials on $U $. The
set of complex valued points $U(\C) $ of $U $ carries a canonical
structure of a smooth analytic manifold which we denote by $\Uan
$. There is a natural morphism from the algebraic de Rham cohomology
to its transcendental counterpart $\Hdr^p(U) \otimes_k \C \to
\Hdr^p(\Uan) $ which is an isomorphism due to a theorem of
A. Grothendieck (\cite{groth}). More generally, let $\nabla:E \to E \otimes
\Omega^1_{U|k} $ denote an integrable algebraic connection on $U $ on
the vector bundle $E $. The integrability means that the associated
sequence
$$
0 \to E \stackrel{\nabla}{\lto} E \otimes_{\O_U} \Omega^1_{U|k} \lto
\ldots \stackrel{\nabla}{\lto} E \otimes_{\O_U} \Omega^n_{U|k} \to 0
$$
is indeed a complex. Again,
there is a natural morphism $\Hdr^p(U; E, \nabla) \otimes_k \C \to
\Hdr^p(\Uan; E^{\rm an},\nabla^{\rm an}) $ between the algebraic and
analytic de Rham cohomology, both defined as the hypercohomology of
the corresponding de Rham complexes. These morphisms fail to be
isomorphisms in general.

In \cite{deligne}, P. Deligne introduces the notion of regular
singularities of the connection at infinity and proves that the
comparison morphism
is an isomorphism in case of a {\em regular} singular
connection. Later Z. Mebkhout gives an alternative proof by
introducing the irregularity sheaf of the connection which
contrary to Deligne's proof does not use Hironaka's resolution of singularities
(\cite{mebkhout}). 

Applying the
analytic Poincar\'e lemma yields another way of stating Deligne's
comparison theorem as the perfectness of the period pairing
\begin{equation}\label{eq:introreg}
(\Hdr^p(U; E, \nabla) \otimes_k \C) \otimes_\C H_p(\Uan, \bve^\wee) \to \C
\end{equation}
induced from integration of differential forms over smooth topological
chains. Here, $H_p(\Uan, \bve^\wee) $ denotes the singular homology of
the analytic manifold $\Uan $ with values in the local system
$\bve^\wee := \ker( E^{{\rm an}, \wee} \to E^{{\rm an}, \wee} \otimes
\Omega^1_{\Xan| \C}) $ of flat sections in the dual bundle $E^{{\rm
    an}, \wee} $ which carries a natural connection $\nabla^\wee $
dual to the given connection on $E $.

Assume that the local system $\bve \subset E^{\rm an} $ of analytic
solutions of $\nabla^{\rm an} $, a locally constant sheaf of $\C
$-vector spaces, comes equipped with the structure of an $F $-local
system for some subfield $F \subset \C $. Then, the singular homology
with values in $\bve $ also carries a natural $F $-structure and the
period pairing gives a well-defined invariant, the alternating product
of the determinants of the pairing for all $p $, as an element in
$k^\times \backslash \C^\times / F^\times $. These invariants for
regular singular connections together with their $\ell $-adic analogs
for tamely ramified sheaves over varieties over finite fields were
extensively studied by T.~Saito and T.~Terasoma in \cite{saitotera}.

In their fundamental paper \cite{b-e}, S.~Bloch and H.~Esnault
generalize the period pairing \eqref{eq:introreg} to the case of {\em
  irregular} singular connections {\em on curves}. To this end, they
define modified homology groups $\Hrd_\ast(\Xan; E^{\rm an},
\nabla^{\rm an}) $ on the
associated Riemann surface which pairs with the algebraic de Rham
cohomology in terms of period integrals and prove perfectness of the
resulting pairing:
$$
(\Hdr^\ast(U; E, \nabla) \otimes_k \C) \otimes \Hrd_\ast(\Xan;
E^{{\rm an}, \wee}, \nabla^{{\rm an}, \wee}) \to \C \ .
$$
The resulting periods are interesting objects by
themselves (the integral representations of the classical
Bessel-functions, Gamma-function and confluent hypergeometric
functions arise in this way as periods of irregular singular
connections on curves) and are mysteriously related to ramification
data for wildly ramified $\ell $-adic sheaves on curves over a
finite field (see e.g. \cite{terasoma}).

In the present paper, we want to start the investigation of the
higher-dimensio\-nal case by studying the period pairing for irregular
singular connections on smooth algebraic surfaces $U $ over $k \subset
\C $. It turns out that additional features arise that do not appear
in the one-dimensional case. These new phenomena affect our proceeding
essentially. We will comment on these after stating the main result.

Let $\nabla^\wee $ be the dual connection. We define the complex
$\Crdt(\nabla^\wee) $ of sheaves of rapid decay
chains on the real oriented blow-up $\widetilde{X} $ of a {\em good}
compactification $(X,D) $ of $\Uan $ (assuming a conjecture of C.~Sabbah, see
section \ref{sec:andrsab} for details), the hypercohomology of which
gives the {\bf rapid decay homology}
$$
\Hrd_p(\Uan; E^\wee, \nabla^\wee) := \Hy^{-p} (\widetilde{X},
\Crdt(\nabla^\wee)) \ .
$$
We then define a natural pairing
\begin{equation} \label{eq:intropair}
\big( \Hdr^p(U; E, \nabla) \otimes_k \C \big) \otimes \Hrd_p(\Uan;
E^\wee, \nabla^\wee) \lto \C
\end{equation}
in terms of period integrals. Our main result, Theorem \ref{thm:main},
asserts that this period pairing is a perfect duality (assuming that
Sabbah's Conjecture holds for $(E, \nabla) $, which is know e.g. if
${\rm rank} \, E \le 5 $).

We are now going to explain the new phenomena arising in the
two-di\-men\-sio\-nal situation as well as Sabbah's Conjecture. For
the situation on curves, the Theorem of
Levelt-Turrittin (cp.~\cite{mal1}) asserts that the formal completion
of the given flat meromorphic connection is locally isomorphic (after
a cyclic covering) to a certain elementary model, namely a direct sum
of irregular singular connections on line bundles times regular
singular connections. 

In dimension two, C.~Sabbah started the investigation of the analogous
questions on the formal structure of meromorphic
connections (cp.~\cite{sabbah1}) indicating a subtle additional
feature: If $\Xan $ denotes an arbitrary
compactification of $\Uan $ with normal crossing divisor $D:= \Xan
\smallsetminus \Uan $ as the complement, the desired formal decomposition
can be expected after a finite sequence of point
blow-ups of $\Xan $ only. Additionally, one asks for the elementary model to
fulfill a certain technical condition, in which case it is called a
{\em good} elementary model. In \cite{sabbah1}, C.~Sabbah conjectures
that such a good formal decomposition can be achieved after a finite
sequence of point blow-ups and cyclic ramification along the smooth
strata of the divisor. He proves the conjecture for several classes
of connections, in particular it holds if the rank of the connection
is lower than or equal to 5.

Assuming Sabbah's Conjecture, we prove a local duality statement
in the following sense. If $X $ is again an
arbitrary compactification of $U $ with a normal crossing divisor $D:=X
\smallsetminus U $, one considers the {\em
  asymptotically flat de Rham complex} $\DRD(\nabla^\wee) $ on the
real oriented blow-up $\widetilde{X} $ of $\Xan $. If $\DRmod(\nabla) $ 
denotes the {\rm de Rham complex with moderate growth} on
$\widetilde{X} $ and $\widetilde{\jmath}:\Uan \hookrightarrow
\widetilde{X} $ the inclusion, there is a natural local duality pairing
$$
\DRmod(\nabla) \otimes_{\C_{\widetilde{X}}} \DRD(\nabla^\wee) \to
\widetilde{\jmath}_!\C_U \ .
$$
We prove that this is a perfect duality (in the derived sense), in
Theorem \ref{thm:locdual}. Taking global sections yields a
global duality result, Theorem \ref{thm:abstrglob}, by
Poincar\'e-Verdier duality which however lacks an explicit description
in terms of periods.

In case of a {\em good} compactification, we prove that there is a
canonical isomorphism in the bounded derived category
\begin{equation}\label{eq:introcrddrd}
\Crdt(\nabla^\wee) \cong \DRD(\nabla^\wee)[2d] \in
D^b(\C_{\widetilde{X}}) \ ,
\end{equation}
where $d:= \dim_\C(X)=2 $. This isomorphism is shown to be compatible
with both the period pairing and Poincar\'e-Verdier duality which
allows to deduce our main result, Theorem \ref{thm:main}, asserting
the perfectness of the period pairing \eqref{eq:intropair}.

Note, that one could try to define the rapid decay sheaves literally
in the same way as in the {\em good} case also for {\em non-good}
compactifications. It seems likely, however, that these sheaves do not
carry enough information to detect the algebraic de Rham cohomology,
since the latter is dual to the hypercohomology of $\DRD(\nabla^\wee)
$ by Theorem \ref{thm:abstrglob}, whereas the proof of the
isomorphism between $\DRD(\nabla^\wee) $ and the rapid decay
complex $\Crdt(\nabla^\wee) $ in the derived category uses the fact
that we start with a good compactification in an essential way. As for
the rapid decay sheaves, it is not clear how they behave under point
blow-ups (in contrast to the asymptotically flat de Rham complex, see
Lemma \ref{lemma:blow}). 

Finally, the explicit description of the rapid decay chains allows us
to forward a given rational structure on the local system to the rapid
decay homology vector spaces. To be more precise, assume that the $\C
$-local system $\bve $ comes equipped with a structure of a local
system of $F $-vector spaces for some subfield $F \subset \C $. The
rapid decay homology then naturally inherits an $F $-structure so that
there is a well-defined determinant of the period pairing
as an element in $k^\times \backslash \C^\times / F^\times $ for
{\em irregular} singular connections on surfaces (see Definition
\ref{def:det}), a generalization in dimension two of T.~Saito and
T.~Terasoma's definition for {\em regular} singular connections (which works
in any dimension).

{\bf Notational convention.} Distinguishing between algebraic varieties and
associated complex manifolds, we will decorate the latter with a
superscript 'an', so that $\Xan $ denotes the complex manifold
associated to the algebraic variety $X $ over $k \subset \C $. For
readability reasons, however, we will not always pursue this notation in the
case of the vector bundles and connections considered whenever it is
clear, e.g. from the type of the space over which they live, in which
category we are working.

\section{Rapid decay homology and periods} \label{sec:rmrk}

\subsection{The geometric situation}

Let $k \subset \C $ be a subfield of the field of complex numbers and
let $U $ be a smooth quasi-projective algebraic surface defined over $k
$. We further consider a vector bundle $E $, i.e. a locally free $\O_X
$-module of rank $r $, together with a flat connection
$$
\nabla: E \to E \otimes_{\O_U} \Omega^1_{U|k}
$$
on $E $. We remark that we do not
impose any condition on the behavior of $\nabla $ at infinity in some
compactification of $U $.

Extending scalars from $k $ to $\C $, we can switch to the analytic
topology. The algebraic connection $\nabla $ induces a flat analytic connection
on the vector bundle $E^{\rm an} $ on the analytic manifold $U^{\rm an} $. We
denote by
$$
\bve := \ker(E^{\rm an} \stackrel{\nabla}{\to} E^{\rm an} \otimes
\Omega_{U^{\rm an}}^1) \subset E^{\rm an} 
$$
the corresponding local system of horizontal sections. By Cauchy's
theorem on linear differential equations for complex variables, the
subsheaf $\bve $ of $E^{\rm an} $ is a locally constant sheaf of $\C
$-vector spaces of the same rank as $E^{\rm an} $.

Let $F \subset \C $ denote another subfield of $\C $. Let us assume
for a moment that
the local system $\bve $ on $U^{an} $ comes equipped with a given $F
$-structure. More precisely, in analogy to \cite{saitotera} we
consider the category $W_{k,F}(U) $ of triples $\M=((E,\nabla),\bve_F,
\rho) $ with: \\[-1.5\baselineskip]
\begin{enumerate} 
\item a vector bundle $E $ on $U $ with rank $r $ together with a flat
  connection $\nabla:E \to E \otimes_{\O_U} \Omega^1_U $, \\[-1.5\baselineskip]
\item a local system $\bve_F $ of $F $-vector spaces on the analytic
  manifold $U^{an} $,\\[-1.5\baselineskip]
\item a morphism $\rho: \bve_F \to E^{an} $ of sheaves on $U^{an} $
  inducing an isomorphism $\bve_F \otimes_F \C \stackrel{\sim}{\to}
  \ker(\nabla^{an}) $ of local systems of $\C $-vector spaces on
  $U^{an} $.
\end{enumerate}
A morphism between $((E, \nabla), \bve_F,
\rho) $ and $((E', \nabla'), {\bve_F}', \rho') $ is given by a
morphism $E \to E' $ respecting the connections together with a
morphism $\bve_F \to {\bve_F}' $ of $F $-local systems with the
natural compatibility condition with respect to $\rho $ and $\rho' $.

In \cite{saitotera}, T.~Saito and T.~Terasoma consider a similar
situation, for arbitrary dimension $\dim(X) $ with the restriction
that $\nabla $ has to be regular singular at infinity (i.e. along the
boundary of a suitably chosen compactification $X $ of $U $). The
perfectness of the resulting period pairing then follows from
P. Deligne's fundamental comparison theorem already mentioned in the
introduction before.

\subsection{Good formal structure (after C.~Sabbah)}
\label{sec:andrsab}

Applying the desingularization theorem of Hironaka, we may assume that
$U $ is embedded in a smooth projective variety $X $ such that $D:=X
\smallsetminus U $ is a divisor with normal crossings. Passing
to the analytic topology, we consider the meromorphic connection
$$
\nabla^\an: E^\an(\ast D) \to E^\an(\ast D) \otimes \Omega^1_{\Xan|\C}
\ .
$$
According to a conjecture of C.~Sabbah, any such connection admits a
good formal structure after a finite sequence of point blow-ups in the
following sense:

First, we recall the definition of a regular singular
connection. Consider the local situation at a point $x^0 \in D
$. Choosing suitable coordinates with $x^0=0 $, we have either $D= \{
x_1x_2=0 \} $ or $D= \{ x_1=0 \} $. Suppose, that $x^0 $ is a
crossing-point of $D $. Then, an integrable $\O_X(\ast D)
$-connection $(R, \nabla) $ is called {\bf regular singular}, if there
is a finite dimensional $\C $-vector space $V $ together with
two commuting endomorphisms $\delta_i:V \to V $, $i=1,2 $, such that
$R $ is locally isomorphic to $\O_{\Xan}(\ast D) \otimes_\C V $ with
$$
\nabla (f \otimes v) = df \otimes v + f \cdot d\log \, x_1 \otimes
\delta_1(v) + f \cdot d\log \, x_2 \otimes \delta_2(v) \ .
$$
In the case of a smooth point of $D $, one asks for the
endomorphism $\delta_1 $ only and requires that $\nabla (f \otimes v)
= df \otimes v + f \, d\log \, x_1 \otimes \delta_1(v) $. In
particular, a regular singular line bundle is locally isomorphic to
the connection on the trivial bundle $\O_{\Xan}(\ast D) $ given by
$$
\nabla_\lambda 1 = \lambda_1 \, d\log \, x_1 + \lambda_2 \, d\log \, x_2
$$
with $\lambda_i \in \C $ and $\lambda_2=0 $ if $x^0 $ is a smooth
point of $D $. Such a connection will be denoted by 
\begin{equation} \label{eq:xlambda}
x^\lambda := (\O_{\Xan}(\ast D), \nabla_\lambda) \ .
\end{equation}
We remark, that it follows directly from the definition that any
regular singular connection is locally isomorphic to a successive
extension of regular singular line bundles.

Now, in dimension one, the classical Levelt-Turrittin Theorem asserts that
any meromorphic connection on a curve formally decomposes into the
direct sum of irregular singular line bundles times regular singular
connections (cp. \cite{mal1}). Sabbah's Conjecture gives an 
analogous formal decomposition property for surfaces with the additional difficulty
that one has to allow preceding blow-ups of points on $D $.

If $\alpha $ is a section in $\O_{\Xan}(\ast D) $, we denote by $e^\alpha $
the meromorphic connection on the trivial line bundle $\O_{\Xan}(\ast D) $
given by $\nabla 1 = d \alpha $ (such that the horizontal sections are
multiples of $e^\alpha $). Its isomorphism class only depends on the
class $\alpha \in \O_{\Xan}(\ast D)/ \O_{\Xan} $. An {\bf elementary
  local model} is a meromorphic connection with poles along $D $
isomorphic to a direct sum
$$
\bigoplus_{\alpha \in A} e^{\alpha} \otimes R_\alpha \ ,
$$
where $A $ is a finite subset of $\O_{\Xan}(\ast D) / \O_{\Xan} $ and the
$R_\alpha $ are regular singular
connections. The model is called {\bf good}, if
\begin{enumerate}
\item the various $\alpha $ are pairwise different in
  $\O_{\Xan}(\ast D) / \O_{\Xan} $. In particular there is at most one
  trivial $\alpha=0 $.
\item For $\alpha \neq \beta \in A $, the divisor $(\alpha-
  \beta) $ has support on $D $ and is negative, and the same holds for
  the divisor $(\alpha) $ for any non-trivial $\alpha \in A $.
\end{enumerate}
Now, for any stratum $Y $ of $D $ in the natural stratification (i.e.
either $Y $ is a connected component of the smooth part of $D $ or a
crossing point), let $\OXY$ be the completion of $\O_{\Xan} $ along $Y $. A
connection $\MM $ meromorphic along $D $ is said to have a {\bf good
  formal decomposition along} $(D,Y) $ {\bf at a point} $x_0 \in Y $,
if locally at $x_0 $ there exists a good elementary model $\MM^{el} $
and an isomorphism
$$
\MM \otimes_{\O_{\Xan}} \OXY \cong \MM^{el} \otimes_{\O_{\Xan}} \OXY
\mbox{\qquad locally at } x_0 \ .
$$

The connection $\MM $ is said to have a {\bf good formal structure
  along} $(D,Y) $ at $x_0 $, if after a bicyclic ramification
along the components of $D $ the inverse image has a good formal
decomposition. Finally, $\MM $ has a {\bf good formal structure along}
$D $, if it has such a structure along $(D,Y) $ for any stratum $Y $
at any point. The precise formulation of Sabbah's Conjecture now reads
as follows.
\begin{conjecture}[C.~Sabbah, \cite{sabbah1} I.2.5.1] \label{conj:sab}
There is a finite sequence of point blow-ups $b:X' \to X $ over $x_0
$, such that the inverse image connection $b^\ast \MM $ has a good formal 
  structure in a neighborhood of $x_0 $.
\end{conjecture}
C.~Sabbah proves his conjecture for several classes of
connections. In particular he achieves the following result.
\begin{theorem}[C.~Sabbah, \cite{sabbah1} I.2.5.2] \label{thm:sab}
The conjecture is true for any connection $\M $ with ${\rm rank} \, \M
\le 5 $.
\end{theorem}
In the following, we will always assume that Sabbah's Conjecture
holds.

\subsection{Definition of rapid decay homology} \label{sec:defofrd}

Recall that we start with an integrable algebraic connection
$\nabla:E \to E \otimes \Omega^1_U $ over the smooth quasi-projective
two-dimensional variety defined over a subfield $k \subset \C
$. Assuming Sabbah's Conjecture, we can
choose an embedding of $U $ into a smooth projective
variety $X $, such that $D:=X \smallsetminus U $ is a divisor with
normal crossings and the associated meromorphic connection on the
complex analytic manifold $X^\an $ admits a good formal structure
locally at any point of $D $. 

Since $D $ is a normal crossing divisor, we can consider the real
oriented blow-up of the irreducible components of $D $. We denote by
$\pi: \widetilde{X} \to \Xan $ the composition of all these. Locally $\pi
$ reads as
$$
\pi: (\R^+ \times S^1)^2 \to \Xan \ , \ (r_1, \vt_1, r_2, \vt_2) \mapsto
(r_1 e^{i \vt_1}, r_2 e^{i \vt_2}) 
$$
in the local case $D=\{ x_1x_2=0 \} $ and
$$
\pi: (\R^+ \times S^1) \times Y   \to \Xan \ , \ (r_1, \vt_1, x_2) \mapsto
(r_1 e^{i \vt_1}, x_2) 
$$
in the case $D=\{ x_1=0 \} $ where $Y $ is an open disc in $\C $. 

Let $\widetilde{D}:= \pi^{-1}(D) $. Then $\pi $ defines a
diffeomorphism between $\widetilde{X} \smallsetminus \widetilde{D} $
with $\Uan $ which we will identify in the following and write
$\widetilde{\jmath}: \Uan \hookrightarrow \widetilde{X} $ for the inclusion.

We consider singular homology of the analytic manifolds
involved. Since we want to study period integrals, we will a priori
work with smooth singular chains (which give the same homology as the
purely topological chains by smooth approximation). In this sense, let
$S_p(\widetilde{X}) $ denote the $\Q $-vector space of
smooth singular $p $-chains, i.e. the free vector space over all piecewise
smooth maps $c:\Delta^p \to \widetilde{X} $ from the standard $p $-simplex to
$\widetilde{X} $. The natural boundary operator $\partial:
S_p(\widetilde{X}) \to S_{p-1}(\widetilde{X}) $ defines a chain
complex the homology of which gives the singular homology. If $A \subset
\widetilde{X} $ denotes a closed subset, the relative chain complex is
defined to be the quotient $S_{\bdot}(\widetilde{X},A) :=
S_{\bdot}(\widetilde{X}) / S_{\bdot}(A) $.

These notions sheafify in the following sense (cp. \cite{verdier},
\cite{swan}). Let $\Ct^{-p} $ denote the sheaf associated to the
presheaf $V \mapsto S_p(\widetilde{X}, \widetilde{X} \smallsetminus V)
$ of $\Q $-vector spaces. The usual boundary operator makes
$\Ct^{-\bdot} $ into a complex. Note that we will use the standard
sign convention, i.e. if we write $\partial $ for the topological boundary operator
on chains, the differential of $\Ct^{-\bdot} $ will be given
by $(-1)^r \partial $ on $\Ct^{-r} $.

Let $d:= \dim_\C(X)=2 $ denote the complex dimension of $X $. Then the
sheaf $\widetilde{\jmath}_! \C_{\Uan}[2d] $ is the dualizing sheaf on the
compact real manifold $\widetilde{X} $ with boundary $\widetilde{X}
\smallsetminus \Uan $ and $\Ct^{-\bdot} $ is a resolution
of this sheaf (cp.~\cite{verdier}).

Let $\Ctd^{-\bdot} $ be the complex of relative chains, i.e. the
sheaf associated to $V \mapsto S_{\bdot}(\widetilde{X}, (\widetilde{X}
\smallsetminus V) \cup \widetilde{D}) $. We are going to define the
complex of rapid decay chains as a 
subquotient of the complex $\Ctd^{-\bdot} \otimes_{\Q}
\widetilde{\jmath}_\ast \bve $ of the complex of sheaves of
relative chains with values in the local system
$\widetilde{\jmath}_\ast \bve $. Let $V \subset \widetilde{X} $
be an open subset with $V \cap \widetilde{D} \neq 0 $ and $c \otimes
\ve $ be a local section of this complex over $V $. Take $c $ to be 
represented by a piecewise smooth map from the standard $p $-simplex
$\Delta^p $ to $\widetilde{X} $. We introduce the notion of rapid decay chains
which works in any dimension $\dim(X)=d $.
\begin{definition} \label{def:rdchains}
We say that the local section $c \otimes \ve $ is a {\bf rapid decay
  $p $-chain} if the section $\ve \in \Gamma(V \smallsetminus
  \widetilde{D}, \bve) $ is {\bf rapidly decaying along $c $ inside $V $} in the
  following sense: For any $y \in c(\Delta^p) \cap \widetilde{D} \cap
  V $, let $e=(e_1, \ldots, e_r) $ denote a local trivialization of
  $E(\ast D) $ and $z_1, \ldots, z_d $ local coordinates such that
  locally $D = \{ z_1 \cdots z_k = 0 \} $ and that $y=0 $. With respect to the
  trivialization $e $ of $E $, $\ve $ becomes an
  $r $-tuple of analytic functions 
$$
  f_i:= e_i^\ast \ve: (c(\Delta^p) \setmin \widetilde{D}) \cap V \to \C \ , \
  (z_1, \ldots, z_d) \mapsto f_i(z) \ .
  $$
  We require that these functions have rapid decay for the argument
  approaching $\widetilde{D} $, i.e.~that
  for all $N \in\N^k $ there is a $C_N > 0 $ such that
  $$
  |f_i(z)| \le C_N \cdot |z_1|^{N_1} \cdots |z_k|^{N_k}
  $$
  for all $z \in (c(\Delta^p) \setmin \widetilde{D}) \cap V $ with small
  $|z_1|, \ldots, |z_k| $.

The subsheaf of $\Ctd^{-p} \otimes_\Q \widetilde{\jmath}_\ast
\bve $ generated by all rapid decay $p $-chains will be denoted by
$ \Crdtp{-p}(\nabla) $. Together with the usual boundary operator of
chains, these give the {\bf complex of rapid decay chains}
$\Crdt(\nabla):=(\Crdtp{-\bdot}(\nabla), \partial) $.
\end{definition}
We stress that we do not impose any decay condition on pairs $(c, \ve)
$ with $c(\Delta^p) \subset \widetilde{D} $; nevertheless we call
those pairs rapidly decaying as well. Note also that the choice of the
meromorphic trivialization does not effect the notion of rapid decay.

The homology we are going to study is
defined as follows. Recall that we assume Sabbah's Conjecture to hold
for the given connection.
\begin{definition}
Let $(E, \nabla) $ be an integrable connection on the smooth
quasi-projective algebraic surface $U $ over $k \subset \C $. The
{\bf rapid decay homology} of $(E, \nabla) $ is defined to be the
hypercohomology
$$
\Hrd_k(\Uan; E, \nabla) := \Hy^{-k}(\widetilde{X}, \Crdt(\nabla)) \ ,
$$
where $\Xan $ is a good compactification of $\Uan $ with respect to $\nabla
$ in the sense of Sabbah's Conjecture and $\pi:\widetilde{X} \to \Xan $
denotes the oriented real blow-up of the normal crossing divisor $D:= \Xan
\smallsetminus \Uan $.
\end{definition}

Later, we will prove that this definition is independent of the choice
of the good compactification (Proposition \ref{prop:indep}),
justifying the notation.

Note, that the usual barycentric subdivision operator obviously can be
defined on the rapid decay complex as well. Hence $\Crdt(\nabla)
$ is a homotopically fine complex of sheaves (cp.~\cite{swan}, p.~87)
and its hypercohomology can be computed as the cohomology of the
complex of global sections.

\subsection{The pairing and statement of the main result} \label{sec:distri}

We are now going to define the pairing between de Rham cohomology of
$\nabla $ and the rapid decay homology of the dual connection
$\nabla^\wee $ on the dual bundle $E^\wee $, characterized by the
equality $d \langle e,\vi \rangle = \langle \nabla
e, \vi \rangle + \langle e, \nabla^\wee \vi \rangle $ for local
sections $e $ of $E $ and $\vi $ of $E^\wee $. Let $\bve^\wee \subset
(E^\an)^\wee $ denote the corresponding local system.

We will work on the real oriented blow-up
$\pi: \widetilde{X} \to \Xan $ of $D $ in a good compactification $(X,D)
$ of $U $ with respect to $\nabla $. For later purposes, we want to
describe the algebraic de Rham cohomology of the connection in terms
of naturally defined sheaves on $\widetilde{X} $. To this end, let
$\Amod $ denote the sheaf of functions on $\widetilde{X} $ which are
holomorphic on $\Uan \subset \widetilde{X} $ and of moderate growth along
$\pi^{-1}(D) $. Then $\Amod $ is flat over $\pi^{-1}(\O_X) $ and if we
define the {\bf moderate de Rham complex} of $(E, \nabla) $ to be
$$
\DRmod(\nabla) := \Amod \otimes_{\pi^{-1}(\O_X)} \pi^{-1}(\DR(\nabla))
\ ,
$$
this complex computes the meromorphic de Rham cohomology of $\nabla $
on $\Xan $ and hence the algebraic de Rham cohomology of $\nabla $ on
$U $ (cp.~with \cite{sabbah3}, Lemma 1.3 or \cite{sabbah1}, Corollaire
1.1.8):
\begin{equation} \label{eq:algmod}
\Hdr^k(U; E, \nabla) \cong \Hy^k (\widetilde{X}, \DRmod(\nabla)) \ .
\end{equation}

We want to define the pairing by means of integration of differential
forms over smooth topological chains as a pairing of complexes of sheaves. To
this end, we are going to define the sheaf of distributions with rapid
decay as follows. First, let
$\Db_{\widetilde{X}}^{-s} $ be the usual sheaf of distributions of
degree $s $ on
$\widetilde{X} $, i.e. the local sections for an open $V \subset
\widetilde{X} $ are the continuous linear functionals
$$
\Db_{\widetilde{X}}^{-s}(V) := \Hom_{\rm cont}( \Gamma_c(V, \Etop{s}),
\C)
$$
on the space $\Etop{s} $ of $C^\infty $ differential forms on
$\widetilde{X} $ of
degree $s $ with compact support in $\widetilde{X} $. The topology on
$\Etop{s} $ is defined to be the limit topology where we write
$\Etop{s} $ as the direct limit of all differential forms with support
in some fixed compact $K \subset V $ with their usual
Fr\'echet-topology, the limit taken over all these $K $
(cp. \cite{hoerm}, chapter II).

We now choose local coordinates $x_1,x_2 $ in $X $ such that locally $D=\{
x=0\} $, where $x=x_1 $ in case of a smooth point of $D $ and
$x=x_1x_2 $ in case of a crossing point of $D $. We call a
distribution $\vi \in \Db_{\widetilde{X}}^{-s}(V) $ a {\bf rapid decay
  distribution} if for any compact $K \subset V $ and any element $N
\in \N $ or $\N^2 $ (depending on the local type of $D $) there exist
$m \in \N $ and $C_{K,N} >0 $ such that for any test form $\eta $
with compact support in $K $ the estimate
\begin{equation} \label{eq:rddistrib}
|\vi(\eta)| \le C_{K,N} \sum_i \sup_{|\alpha| \le m} \sup_K \{ |x|^N
|\partial^\alpha f_i| \}
\end{equation}
holds, where $\alpha $ runs over all multi-indices of degree less than
or equal to $m $ and $\partial^\alpha $ denotes the $\alpha $-fold
partial derivative of the coefficient functions $f_i $ of $\eta $ in
the chosen coordinates.

Let $\Dbrd{-s} $ denote the resulting sheaf of rapid decay
distributions on $\widetilde{X} $. Varying $s $, we obtain a
complex of sheaves where once again the sign convention applies,
i.e. the differential reads as
$$
\Dbrd{-s} \lto \Dbrd{-s+1} \ , \ \vi \mapsto \big( \eta \mapsto
(-1)^s \vi(d\eta) \big) \ .
$$
Standard arguments show that the resulting complex $\Dbrd{-\bdot} $ is a fine
resolution of the
extension $\widetilde{\jmath}_! \C_{\Uan}[2d] $ of the constant sheaf by
$0 $ along $\widetilde{\jmath}:\Uan \hookrightarrow \widetilde{X} $
shifted by the real dimension $2d $. 



\medskip
Now, consider a local section $\omega $ of $\DRmod $ in degree $s
$, as well as a rapid decay chain $c \otimes \ve \in \Gamma(V,
\Crdtp{-r}(\nabla^\wee)) $ with respect to the dual connection, $c $
being a smooth topological $r $-simplex $c $ in $\widetilde{X} $. By
definition, $\ve $ is rapidly decaying along $c $ and $\omega $ has at
most moderate growth along $\widetilde{D} $, hence the local section
$\langle \ve, \omega \rangle \in \Gamma(V \smallsetminus
\widetilde{D}, \Omega^{\infty,s}_{\Uan}) $ remains rapidly decaying
along $c $ and the same holds for
$$
\eta \wedge \langle \ve, \omega \rangle  \in \Gamma(V \cap
\Uan, \Omega_{\Uan}^{\infty,r})
$$
for any test form $\eta \in \Gamma_c(V, \Etop{p}) $ with $p=r-s
$. Therefore, the integral
\begin{equation} \label{eq:integralex}
\int_{c} \eta \wedge \langle \ve, \omega \rangle
\end{equation}
obviously converges. The rapid decay condition of $\ve $ also ensures
that the distribution defined by \eqref{eq:integralex} satisfies the
estimate \eqref{eq:rddistrib}, since the coordinate functions of $\ve
$ along the curve of integration $c $ can be bounded from above by any
power of $|x| $ by definition. For support reasons, the integral is
well-defined also, i.e. independent of the choice of $c $ in its 
equivalence class modulo chains in $\widetilde{X} \smallsetminus V $. 

The above considerations define a morphism of sheaves
\begin{equation} \label{eq:integral}
\DRmd{s}(\nabla) \otimes \Crdtp{-r}(\nabla^\wee) \to \Dbrd{s-r} \ ,
\end{equation}
which indeed induces a morphism of complexes. More precisely, assuming
that $c^{-1}(\widetilde{D}) \subset \partial \Delta^p $ (which is no
restriction due to subdivision) if $c_t $
denotes the chain one obtains by cutting off a small
tubular neighborhood with radius $t $ around the boundary $\partial
\Delta^p $ from the given chain $c $. Then, for $c \otimes
\ve $ as before and $\omega $ a form with moderate growth of degree
$(s-1) $ and $\eta $ as before, we have the {\bf 'limit Stokes formula'}
\begin{multline}\label{eq:paircompl}
\Int_c \eta \wedge \langle \ve, \!\! \nabla \omega \rangle + (-1)^{r+s}
\!\! \Int_{\partial c} \eta \wedge \langle \ve, \omega \rangle  = \lim_{t \to 0}
(\Int_{c_t} \eta \wedge \langle \ve, \nabla \omega \rangle  +
(-1)^{r+s} \!\! \Int_{\partial c_t} \eta \wedge \langle \ve, \omega \rangle )
= \\
\stackrel{\text{Stokes}}{=} 
(-1)^{r-s} \lim_{t \to 0} \Int_{c_t} d \eta \wedge \langle \ve, \omega
 \rangle = (-1)^{r-s} \Int_c d\eta \wedge \langle \ve, \omega
 \rangle
\end{multline}
where one has to keep in mind that by the given growth/decay
conditions the integrals over the faces of $\partial c_t $ 'converging'
against the faces of $\partial c $ contained in $\widetilde{D} $
vanish. Note also, that by the definition of the dual connection $d
\langle \ve, \omega \rangle = \langle \ve, \nabla \omega \rangle $,
since $\nabla^\wee \ve = 0 $, and that the sign conventions in
defining the total complex of the distributions, the rapid decay
complex and in the rule for differentiating a wedge product of forms
give the appropriate sign in \eqref{eq:paircompl}.

Equation \eqref{eq:paircompl}, however,
directly shows that \eqref{eq:integral} is compatible with the
differentials of the complexes involved, i.e. we have obtained a
pairing of complexes
$$
\DRmod(\nabla) \otimes_\C
\Crdt(\nabla^\wee) \to \Dbrd{-\bdot} \ .
$$


Since $\Dbrd{-\bdot} $ is a resolution of $\widetilde{\jmath}_!
\C_{\Uan}[2d] $ it follows that $\Hy^0(\widetilde{X}, \Dbrd{-\bdot}) \cong \C $,
keeping in mind that $\widetilde{\jmath}_! \C_{\Uan}[2d] $ is the
dualizing sheaf for the compact real manifold $\widetilde{X} $ with
boundary.

Taking hypercohomology of the above pairing in degree $0 $ thus induces a pairing
\begin{equation} \label{eq:perfV}
  H^p_{dR}(U; E,\nabla) \otimes_\C \Hrd_p(\Uan; (E^\wee, \nabla^\wee))
  \lto \C 
  \end{equation}
will be called the {\bf period pairing} of the algebraic connection
$(E, \nabla) $. Our main result is the following
\begin{theorem}[Global duality of the period pairing] \label{thm:main}
Let $(E, \nabla) $ be a flat connection on the smooth
quasi-projective algebraic surface $U $ over the subfield $k \subset
\C $. Assume that Sabbah's Conjecture holds for $(E, \nabla) $. Then the pairing
\eqref{eq:perfV} 
$$
\big( \Hdr^p(U; E, \nabla) \otimes_k \C \big) \otimes_\C \Hrd_p(\Uan;
E^\wee, \nabla^\wee) \lto \C
$$
is a perfect pairing of finite-dimensional $\C $-vector spaces.
\end{theorem}

\begin{remark}\label{rmrk:int}
\setcounter{remark:int}{\value{definition}}{\bf i)} 
For ${\rm dim}(X)=1 $ the analogous statement was proven by S.~Bloch
and H.~Esnault in \cite{b-e}. Note
that in the one-dimensional situation any compactification will
suffice for the definition of rapid decay homology since the given
integrable connection always admits a formal structure by the
classical Levelt-Turrittin Theorem and this formal structure is
automatically good. In the two-dimensional situation, the suitable
choice of good compactification plays an important role for the
definition of the rapid decay homology. We will come back to this
later (see section \ref{subsec:rdundcomp}).

\noindent{\bf ii)} We used the de Rham complex with moderate
growth on $\widetilde{X} $ in order to define the period pairing
referring to \eqref{eq:algmod}. However, one can obviously replace
this complex by the pull-back $\pi^{-1}(\DR(E(\ast D),\nabla)) $ in the
definition to obtain the same period pairing.

\noindent{\bf iii)} In case $U $ is affine, the period pairing admits a
direct description in terms of period integrals as follows. For affine
$U $, the algebraic de Rham cohomology can be computed by the
cohomology of global sections in $U $, i.e.
$$
\Hdr^p(U; E, \nabla) = H^p( \ldots \to \Gamma_U(E \otimes
\Omega_{U|k}^q) \stackrel{\nabla}{\to} \Gamma_U(E \otimes
\Omega_{U|k}^{q+1}) \to \ldots) \ .
$$
Using the barycentric subdivision operator, the rapid decay
complex is easily seen to be homotopically fine and thus the rapid
decay homology can be computed taking global sections
of the rapid decay complex also. The period pairing then obtains the
following shape
$$
\Hdr^p(E, \nabla) \otimes_\C \Hrd_p(X, \nabla^\wee) \to \C \ , \
([\omega], [c \otimes \ve]) \mapsto \int_c \langle \ve, \omega \rangle
$$
for a flat algebraic $p $-form $\omega $ on $U $ and a rapid decay cycle
$c \otimes \ve $ on $\widetilde{X} $.

\noindent
{\bf iv)} Due to Sabbah's Theorem \ref{thm:sab}, the conclusion of the
above theorem holds unconditionally for ${\rm rank} \, E \le 5 $.
\end{remark}

\subsection{The determinant of periods}

\noindent
Consider an element $((E,\nabla), \bve_F, \rho) \in W_{k,F}(U) $ for
given subfields $k,F \subset \C $. Since $(E, \nabla) $ and $U $ are
defined over $k $, the de Rham cohomology $\Hdr(U,
\nabla) $ is a $k $-vector space by definition.

The $F $-structure $\bve_F $ on the local system of horizontal
sections given by $\rho $, which obviously also defines an $F $-structure
$\bve^\wee_F $ on the local system $\bve^\wee $ of the dual
connection, induces a canonical $F $-structure on the
complex of rapid decay chains on a chosen good compactification $X
$. More precisely, we define
$$
\Crdtp{-p}(\nabla^\wee)_F \subset \Ctd^{-p} \otimes_{\Q}
\widetilde{\jmath}_\ast \bve^\wee_F
$$
to be the sheaf of subvector spaces over $F $ generated by all rapidly
decaying chains in $\Ctd^{-p} \otimes_{\Q} \widetilde{\jmath}_\ast
\bve^\wee_F $, the property of rapid decay literally being the same as in
definition \ref{def:rdchains}. Its
hypercohomology gives the $F $-vector space $\Hrd_p(\Uan, E, \nabla)_F $
such that 
$$
\Hrd_p(\Uan, E, \nabla)_F \otimes_F \C \stackrel{\cong}{\lto} \Hrd_p(\Uan,
E, \nabla) \ ,
$$
the isomorphism induced by $\rho $. In summary, we can define the
determinant of the period pairing as follows:
\begin{definition}\label{def:det}
  For $((E, \nabla), \bve_F, \rho) \in W_{k,F}(U) $, we define its
  {\bf period determinant} to be the element
$$
\det((E,\nabla), \bve_F, \rho) := \prod_{p \ge 0} \det ( \langle
\gamma^{(p)}_j , \omega^{(p)}_i \rangle )_{i,j}^{(-1)^p} \in k^\times
\backslash \C^\times / F^\times \ , 
$$
where $\omega^{(p)}_i $ denotes a basis of $\Hdr^p(U, E, \nabla) $
over $k $ and $\gamma^{(p)}_j $ a basis of the $F $-vector space
$\Hrd_k(\Uan, E^\wee, \nabla^\wee)_F $.  
\end{definition}
Obviously, the determinant does not depend on the choices made. For
regular singular connections $(E, \nabla) $ this definition coincides
with the one in \cite{saitotera}. In
case $U $ is affine, the matrices involved carry actual period
integrals as entries (see Remark \ref{sec:rmrk}.\arabic{remark:int} ii). 

\section{Local duality}

We are going to study a local duality pairing on the real blow-up
$\widetilde{X} $ of a good compactification $\Xan $ of $\Uan $ with respect
to the given flat connection $(E, \nabla) $ on $U $. Our
main result, Theorem \ref{thm:main}, will follow from the local
duality by standard globalization and the comparison of the period
pairing with the local pairing.

\subsection{Sheaves of functions on the real oriented blow-up}

We identify the complex of rapid decay sheaves with the
asymptotically flat de Rham complex on the real oriented blow-up
$\widetilde{X} $ in a given good compactification $X $ of $U $ which
we are going to define in this section. We first recall the definition of
the following sheaves of functions on the real oriented blow-up. We
can assume a more general situation, namely any compactification $X $
such that $D:= X \smallsetminus U $ is normal crossing. Let
$\widetilde{X} $ denote the real oriented blow-up of $D $. For the
proofs of the various properties stated in the following we refer to
\cite{sabbah1}, II.1. 
\begin{enumerate}
\item The logarithmic differential operators as well as their
  conjugates act on the sheaf $C^\infty_{\widetilde{X}} $ of
  $C^\infty $ functions on $\widetilde{X} $ and one defines
  $$
  \A := \ker \widebar{x}_1 \cdot \widebar{\partial}_{x_1} \cap \ker
  \widebar{x}_2 \cdot \widebar{\partial}_{x_2} \subset
  C^\infty_{\widetilde{X}} 
  $$
in the case $D=\{ x_1x_2=0 \} $ and $\A=\ker \widebar{x}_1 \cdot
\widebar{\partial}_{x_1} \cap \ker \widebar{\partial}_{x_2} $ in the
case $D=\{ x_1=0 \} $. 
  The elements of $\A $ are holomorphic functions on $\widetilde{X} $
  which admit an asymptotic development in the spirit of Majima (cp.
  Proposition B.2.1 in \cite{sabbah1} and \cite{majima}).
\item Let ${\cal P}_{\widetilde{X}}^{<D} $ denote the sheaf
  of $C^\infty $-functions on $\widetilde{X} $ which are flat on
  $\pi^{-1}(D) $, i.e. all of whose derivations vanish on $\pi^{-1}(D)
  $ (cp. \cite{mal3}) and let
  $$
  \AD := \A \cap {\cal P}_{\widetilde{X}}^{<D} \ .
  $$
  The elements of $\AD $ are the holomorphic functions with vanishing asymptotic
  development (cp. Proposition II.1.1.11 in \cite{sabbah1}), i.e. which
  are rapidly decaying on any compactum in $\widetilde{X} $. More
  precisely, if $u $ is a local section of $\AD $ defined on some open
  subset $\Omega \subset \widetilde{X} $, then for any compact $K
  \subset \Omega $ and any $N \in \N^2 $, the function $u $ satisfies
  an estimate of the form
\begin{equation}\label{eq:ad=rd}
|u(x)| \le C_{K,N} \cdot |x_1|^{N_1} |x_2|^{N_2}  \mbox{ for all
  } x \in K \smallsetminus \pi^{-1}(D) \ .
\end{equation}

\item Let $\fA $ denote the formal completion of $\A $ along
  $\pi^{-1}(D) $ and $T_D $ the natural morphism $\A
  \stackrel{T_D}{\lto} \fA := \varprojlim_k \A / {\cal I}_D^k \A
  $.
\end{enumerate}
According to Majima, the sequence $0 \to
\AD \to \A \stackrel{T_D}{\to} \fA \lto 0 $ is exact, generalizing the
analogous theorem of Borel-Ritt in dimension one. 

For $\dim(X)=2 $, both sheaves $\A $ and $\AD $ are flat
as $\pi^{-1}(\O_X) $-algebras. Additionally, we will make use of the
fact that
\begin{equation} \label{eq:ADresol}
\AD \hookrightarrow \big( {\cal P}_{\widetilde{X}}^{<D}
\otimes_{\pi^{-1} C_X^\infty} \pi^{-1} \Epq{0}{\bdot}, \dolb \big)
\end{equation}
is a resolution of $\AD $, where $\Epq{p}{q} $ denotes the sheaf of
$C^\infty $-forms of degree $(p,q) $ on $X$ (cp. Lemme II.1.1.18 in
\cite{sabbah1}). 

Recall that $\Amod $ denotes the sheaf of functions on $\widetilde{X}
$ which are holomorphic on $\Uan $ and of moderate growth along
$\widetilde{D} $. If $\Emd $ denotes the sheaf of $C^\infty
$-functions on $\widetilde{X} $ with moderate growth at $\pi^{-1}(D)
$, the inclusion defines a resolution
\begin{equation} \label{eq:Amodresol}
\Amod \hookrightarrow \big( \Emd \otimes_{\pi^{-1} C^\infty_X}
\pi^{-1} \Epq{0}{\bdot} , \dolb \big) \ .
\end{equation}
Both resolutions being constructed with $C^\infty $-functions consist
of fine sheaves. 

\subsection{Formal classification and asymptotic developments}
\label{sec:sabbahsec}

Given a good formal structure of $(E, \nabla) $ on $X $, the resulting
formal decomposition can be lifted to an asymptotic one in the following
sense. Consider the local situation at a crossing-point and assume $X $
is a small bi-disc around the crossing point $0 $ with coordinates $x_1,
x_2 $ such that locally $D= \{ x_1 x_2= 0 \} $. The situation
at a smooth point is similar with a few obvious changes.

For an $\O_X(\ast D) $-connection $\MM $, let $\MM_{\widetilde{X}} :=
\A \otimes_{\pi^{-1} \O_X} \pi^{-1} \MM $. In the same situation as
above, one says that $\MM $ has a {\bf good ${\cal A} $-decomposition
  along $(D,Y) $ at} $x_0 $ if there exists a good elementary model
$\MM^{el} $ in a neighborhood of $x_0 $ and for all $\vt \in
\pi^{-1}(x_0) $ an isomorphism of the stalks
$$
\MM_{\widetilde{X}, \vt} \cong \MM^{el}_{\widetilde{X}, \vt} \ ,
$$
such that the induced formal isomorphism $\widehat{\MM} \cong
\widehat{\MM^{el}} $ is independent of $\vt $, where one has to keep
in mind that $\fA|_{\pi^{-1}D} = \pi^{-1}
\O_{\widehat{X|D}}|_{\pi^{-1}D} $.  The notions of a {\bf good ${\cal
    A} $-structure} is defined in an analogous manner as above. One
has the following result:
\begin{theorem}[C.~Sabbah, \cite{sabbah1} II.2.1.1]
  \label{thm:sabbahasym} If $\MM $ has a good
  formal decomposition along $(D,Y) $ at $x_0 $, then it lifts to a
  good ${\cal A} $-decomposition at $x_0 $.
\end{theorem}
We remark, that it is essential in the proof of this theorem that the
given formal decomposition is good. Hence, even if $(E, \nabla) $
admits a formal decomposition at $x_0 $ without any preceding blow-up,
it may be necessary to insert point blow-ups in order to arrive at a
good formal decomposition and then be able to lift it to an
asymptotic decomposition. 

\subsection{The asymptotically flat de Rham complex}

\begin{definition}
 Let $\widetilde{X} $ be the real oriented blow-up of a good
 compactification $\Xan $ of $\Uan $ with respect to the given flat
 connection $(E, \nabla) $ on $U $. The complex
$$
\DRD(\nabla^\wee) := \AD \otimes_{\pi^{-1}(\O_{\Xan})}
\pi^{-1}(\DR(\nabla^\wee)) \in D^b(\C_{\widetilde{X}}) 
$$
will be called the {\bf asymptotically flat de Rham complex} of
$\nabla^\wee $.
\end{definition}
Recall that the moderate de Rham $\DRmod(\nabla) $ complex was defined
in a similar way to be
$$
\DRmod(\nabla) := \Amod \otimes_{\pi^{-1}(\O_{\Xan})} \pi^{-1}(\DR(\nabla))
\ .
$$
Now, if $(X,D) $ is good with respect to $(E, \nabla) $, these
complexes simplify as follows:
\begin{proposition} \label{prop:deg0only}
  If $(E, \nabla) $ has a good formal structure along $D $, both
  complexes $\DRmod(E, \nabla) $ and $\DRD(E^\wee, \nabla^\wee) $ have
  cohomology in degree $0 $ only, i.e. the inclusions
\begin{align*}
\Smod & := \H^0( \DRmod(E, \nabla) ) \simeq \DRmod(E, \nabla)  \\
\SD & := \H^0( \DRD(E^\wee, \nabla^\wee)) \simeq \DRD(E^\wee,
\nabla^\wee)
\end{align*}
are quasi-isomorphisms.
\end{proposition}
\pf In \cite{sabbah2}, \S 7, the assertion about $\DRD $ is stated
assuming a {\em very good} formal structure (see the Remark after
Th\'eor\`eme 7.2 and 7.3 in \cite{sabbah2}). It follows from
Th\'eor\`eme II.2.1.2 in \cite{sabbah1}, however, that the same proof
holds if one has a good formal structure only. The proof relies on an
existence theorem for flat solutions to a certain type of partial
differential equations whose entries are rapidly decaying as well. In
the case of $\DRmod $ we will analogously reduce to a theorem on
the existence of solutions with moderate growth, assuming the entries
in the differential equations have moderate growth. This existence
theorem will then be proven in the appendix. 

Since $(E, \nabla) $ is assumed to have a good formal
structure, there exists a bicyclic ramification $\rho:Y \to X $
such that locally on $\widetilde{Y} $ the pull-back connection
$\rho^{-1}(\nabla) $ is isomorphic to its elementary model. Let
$\pi_X:\widetilde{X} \to \Xan $ and $\pi_Y:\widetilde{Y} \to Y^{\rm an} $ denote
the oriented real blow-up of $S:=\rho^{-1}(D) $ and $D $ respectively. Lifting
$\rho $ to $\widetilde{\rho}:\widetilde{Y}
\to \widetilde{X} $, the projection formula yields
\begin{equation} \label{eq:rhoDR}
\bR \widetilde{\rho}_\ast \DRmY (\rho^{-1} \nabla) = \bR
\widetilde{\rho}_\ast \AmdY \Lotimes_{\pi^{-1}(\O_X)} \pi_X^{-1}
\DR(\nabla) \ .
\end{equation}
Now, $\widetilde{\rho} $ being a finite map and since obviously $\bR
\widetilde{\rho}_\ast \AmdY= \Amod $ (using the resolution
\eqref{eq:Amodresol}), it follows that it suffices to
prove the claim on $\widetilde{Y} $. Hence, we can assume 
that $(E, \nabla) $ itself decomposes locally on
$\widetilde{X} $ as
$$
\pi^{-1}(\nabla) \cong_{\rm loc.} \bigoplus\nolimits_{\alpha \in A} e^\alpha
\otimes R_\alpha \ .
$$
Our claim is local in nature, so that we can assume that
$\pi^{-1}(\nabla) = e^\alpha \otimes R_\alpha $.

We consider the local situation on $\widetilde{X} $ above a point $x^0
\in D $. Assume first that $x^0=0 $ is a crossing point, i.e. locally on
$\Xan $ we have the situation $D=\{ x_1x_2=0 \} $. Let $\vt^0 \in
\pi^{-1}(0) $ be the direction under consideration.

Since every regular singular connection is a successive extension of
regular singular line bundles, we can further reduce to the case
$R_\alpha=x^\lambda $ with a $\lambda \in \C^2 $. Then
\begin{equation} \label{eq:drmodvt}
\DRmvt{\vt^0}(\nabla) \cong \Big( 0 \to \Amodvt{\vt^0} \stackrel{\big(
  \begin{smallmatrix} P_1\\ P_2 \end{smallmatrix} \big)}{\lto}
(\Amodvt{\vt^0})^2 \stackrel{(P_2, -P_1)}{\lto} \Amodvt{\vt^0} \to 0
\Big) \ ,
\end{equation}
where
$$
P_i u := x_i \frac{\partial}{\partial x_i} u + x_i \frac{\partial
  \alpha}{\partial x_i} \cdot u + \lambda_i u \ .
$$
The vanishing of $\H^2 $ amounts to solving the partial differential
equation
\begin{equation} \label{eq:de}
(\Sigma_1): \qquad x_1 \frac{\partial}{\partial x_1} u = -x_1 \frac{\partial
  \alpha}{\partial x_1} \cdot u - \lambda_1 u + \rho_1 
\end{equation}
for any given $\rho_1 \in \Amodvt{\vt^0} $, as then $(0,u) $ gives a
preimage of $\rho_1 $. The vanishing of
$\H^1 $ additionally asks for a solution $u $ of \eqref{eq:de} which
also solves
$$
(\Sigma_2): \qquad x_2 \frac{\partial}{\partial x_2} u = -x_2 \frac{\partial
  \alpha}{\partial x_2} \cdot u - \lambda_2 u + \rho_2 \ ,
$$
where $\rho_2 $ is another element in $\Amodvt{\vt^0} $ such
that the system is {\bf integrable} in the sense that
\begin{equation} \label{eq:sysintegrble}
P_2 \rho_1 = P_1 \rho_2
\end{equation}
holds. By Theorem \ref{thm:appendix}, which is
proven in the appendix, such solutions $u \in \Amodvt{\vt^0} $ always
exist. 

In the local case $D=\{ x_1=0 \} $ at $x^0=0 $, we similarly may
assume that $R_\alpha =x_1^{\lambda_1} $ with $\lambda_1 \in \C $ and
then \eqref{eq:drmodvt} remains valid with the same definition for
$P_1 $ and now
$$
P'_2 u:= \frac{\partial u}{\partial x_2} + \frac{\partial
  \alpha}{\partial x_2} \cdot u \ .
$$
The system of partial differential equations to be solved in
$\Amodvt{\vt^0} $ is therefore given by the same equation $(\Sigma_1)
$ as above and 
$$
(\Sigma'): \qquad \frac{\partial}{\partial x_2} u = - \frac{\partial
  \alpha}{\partial x_2} \cdot u + \rho_2 \ ,
$$
where the integrability condition now reads as $P_1 \rho_2= P'_2
\rho_1 $. Again, Theorem \ref{thm:appendix}
proves the existence of such a solution.\\
\qed

\subsection{Local duality for good compactifications} 

Again, let $(E, \nabla) $ be an integrable connection on
$U $ and $(X,D) $ a compactification of $U $ with normal crossing
divisor $D $. Let $\pi:\widetilde{X} \to \Xan $
be the oriented real blow-up of $X $ along $D $. Multiplying an
element of $\AD $, i.e. a function on $\widetilde{X} $ with rapid
decay at $\pi^{-1}(D) $, with a function with moderate growth, results
in a rapidly decaying function again, i.e. multiplication defines a
map $\Amod \otimes_\C \AD \to \AD $. As a consequence, given a local
section
$$
\omega \in \DRmd{p}(E, \nabla) = \Amod \otimes_{\pi^{-1}(\O_{\Xan})}
\pi^{-1} \Omega_{\Xan}^p \otimes_{\pi^{-1}(\O_{\Xan})} \pi^{-1}E \ ,
$$
i.e. a $p $-form with moderate growth, and a local section
$$
\eta \in \DRd{q}(E^\wee, \nabla^\wee) = \AD \otimes_{\pi^{-1}(\O_{\Xan})}
\pi^{-1} \Omega_{\Xan}^q \otimes_{\pi^{-1}(\O_{\Xan})} \pi^{-1}E^\wee \ ,
$$
i.e. a $q $-form with rapid decay, the wedge product $\omega \wedge
\eta $ is rapidly decaying as well. This leads to the following
definition.
\begin{definition}
  The natural pairing
  \begin{equation} \label{eq:locdu} \DRmod(E, \nabla) \otimes_\C
    \DRD(E^\wee, \nabla^\wee) \to \DRD(\O_X, d)
  \end{equation}
  induced by the wedge product is called the {\bf local
    duality pairing} of $(E, \nabla) $.
\end{definition}

Applying Proposition \ref{prop:deg0only} to the trivial connection
$(\O_X,d) $ (for which any compactification is good) gives
$$
\widetilde{\jmath}_! \C_{X
  \smallsetminus D} = \H^0 \DRD(\O_X,d)  \simeq \DRD(\O_X,d) \ ,
$$
with $\widetilde{\jmath}:\Uan \hookrightarrow
\widetilde{X} $ denoting the inclusion.

Now assume that $(X, D) $ is a good compactification with respect to
$(E, \nabla) $. Again by Proposition \ref{prop:deg0only}, the local
duality pairing reduces to the pairing 
\begin{equation}
  \label{eq:locdugood}
  \H^0 \DRmod(E, \nabla) \otimes_\C \H^0 \DRD(E^\wee, \nabla^\wee) \lto
  \widetilde{\jmath}_! \C_{\Uan} \ .
\end{equation}
In this situation, the local duality statement reads as follows:
\begin{theorem}[Local Duality for good formal structure]
  \label{prop:locdugood} \label{thm:locdugood}
If the connection $(E, \nabla) $ admits
  a good formal structure on $X \smallsetminus D $, the pairing
  \eqref{eq:locdugood} 
$$
\Smod \otimes_\C \SD \to \widetilde{\jmath}_! \C_{\Uan}
$$
is a perfect duality, i.e. the induced morphisms
$$
\Smod \to \bR \cHom_{\widetilde{X}}(\SD, \widetilde{\jmath}_!
\C_{X \smallsetminus D}) \mbox{ and } 
\SD \to \bR \cHom_{\widetilde{X}}(\Smod, \widetilde{\jmath}_!
\C_{X \smallsetminus D})
$$
are isomorphisms in $D^b(\C_{\widetilde{X}}) $.
\end{theorem}
\pf Being a local problem on $\widetilde{X} $, we can assume that
after a bicyclic ramification $\rho:Y \to X $ at the point $x_0 \in
D $ we consider, the lift of the pull-back
$\rho^{-1}\nabla $ to the oriented real-blow up $\widetilde{Y} $ of $Y
$ is locally isomorphic to its elementary model. As in the proof of
Proposition \ref{prop:deg0only}, the projection formula applied to
either $\DRD $ or $\DRmod $ again gives
(see \eqref{eq:rhoDR}) $\bR \widetilde{\rho}_\ast \DRmY (\rho^{-1}
\nabla) \cong \DRmod (\nabla) $ and $\bR \widetilde{\rho}_\ast \DRD
(\rho^{-1} \nabla^\wee) \cong \DRD(\nabla^\wee) $, both isomorphism
obviously compatible with the duality morphism given by the wedge
product. With the notation $S:= \rho^{-1}(D) $ and $i: Y
\smallsetminus S \hookrightarrow \widetilde{Y} $ and ${}^\wee \! {\cal
  S}_{\widetilde{Y}}^{<S} := \H^0 \DRdY(\rho^{-1}
\nabla^\wee)) $, we have
$$
\bR \cHom \big( \bR \widetilde{\rho}_\ast {}^\wee \! {\cal
  S}_{\widetilde{Y}}^{<S}, \widetilde{\imath}_! \C_{Y \smallsetminus S} \big) \cong \bR
\widetilde{\rho}_\ast \bR \cHom \big( \SD, \widetilde{\jmath}_! \C_{X
  \smallsetminus D} \big)
$$
by the local Poincar\'e-Verdier duality, where one has to keep in mind
that $\widetilde{\jmath}_! \C $ is up to shift by the real dimension
of $\widetilde{X} $ the dualizing sheaf on the real
manifold $\widetilde{X} $ with boundary $\pi^{-1}(D) $ (including
corners). Therefore it 
remains to prove the claim on $\widetilde{Y} $, i.e. we can assume
that $\nabla $ itself decomposes locally on $\widetilde{X} $. As in
Proposition \ref{prop:deg0only}, we can further reduce to the case
$\pi^{-1}(\nabla) = e^\alpha \otimes x^\lambda $ for some $\alpha \in
\O_{\Xan}(\ast D) / \O_{\Xan} $ and some $\lambda \in \C^2 $.

We distinguish the two case, $x_0 \in D $ being a crossing-point or a
smooth point, i.e. with local coordinates centered at $x^0 $ either
$D=\{ x_1x_2=0 \} $ or $D=\{ x_1=0 \} $. Consider the case of a
crossing-point first.

We define the {\bf Stokes directions} of the elementary connection
$e^\alpha $ with $\alpha(x):= x_1^{-m_1} x_2^{-m_2} u(x) $ with $u(0)
\neq 0 $ (which can be achieved since the decomposition is good) at
the point $0 $ as the elements of 
$$
  \Sigma^0_\alpha := \{ (\vt_1, \vt_2) \mid -m_1 \vt_1 -m_2 \vt_2 + {\rm
    arg}(u(0)) \in (\frac{\pi}{2}, \frac{3 \pi}{2}) \} \subset
  \pi^{-1}(0) \cong S^1 \times S^1 \ .
  $$
Now, consider a given direction $\vt=(\vt_1,\vt_2) \in \pi^{-1}(0)
$. If $\vt \in \Sigma^0_\alpha $, then the solution $x^{-\lambda}
e^{-\alpha(x)} \in \H^0 ( \pi^{-1} \DR(\nabla^\wee)) $ is not rapidly
decaying in any small enough open bisector $V \subset {\widetilde{X}} $ around $\vt
$ whose directions ${\cal V} \subset \pi^{-1}(0) $ are contained in
$\Sigma^0_\alpha $ (recall that $u $ is continuous at $0 $). More
precisely, it is even rapidly growing if $m_1
\neq 0 $ or $m_2 \neq 0 $ and of moderate growth for $m_1=m_2=0 $. Hence 
\begin{equation} \label{eq:sdsigma}
\SD|_V = \widetilde{\jmath}_! \big( x^{-\lambda} e^{-\alpha(x)}
\C_{X \smallsetminus D} \big)|_V
\end{equation}
for any such $V $. Since $x^\lambda e^{\alpha(x)} $ is rapidly
decaying in $V $ for $m_1 \neq 0 $ or $m_2 \neq 0 $ and of moderate
growth otherwise, one has $\Smod|_V = x^\lambda e^{\alpha(x)}
\cdot \C_V $ and consequently the induced morphism
\begin{equation} \label{eq:locV}
\begin{CD}
\bR \cHom ( \SD, \widetilde{\jmath}_!
\C_{X \smallsetminus D} )|_V @>>> \Smod|_V  \\
@V{\cong}VV @VV{\cong}V \\
\bR \cHom ( \widetilde{\jmath}_! \C_{V
  \smallsetminus D} , \widetilde{\jmath}_! \C_{V \smallsetminus D} )|_V
@>{\cong}>> \C_V
\end{CD}
\end{equation}
is indeed an isomorphism. Note hereby, that for the open embedding
$\widetilde{\jmath} $, the functor $\widetilde{\jmath}_! $ is exact and
has $\widetilde{\jmath}^\ast $ as its right adjoint. This easily shows
that the bottom line is an isomorphism.

If $\vt \in \Sigma^0_{-\alpha} $ the situation is similar. If $V $ is a
subsector around $\vt $, whose directions are contained in
$\Sigma^0_{-\alpha} $, then, if $V $ is small enough,
$$
\SD|_V \cong x^{-\lambda} e^{-\alpha(x)} \cdot \C_V
\mbox{\quad and \quad} 
\Smod|_V \cong \widetilde{\jmath}_! \big( x^\lambda e^{\alpha(x)}
\cdot \C_{X \smallsetminus D} \big)|_V \ ,
$$
and \eqref{eq:locV} holds again.

Finally, consider the case that $\vt $ separates the
Stokes directions of $\alpha $ and $-\alpha $, i.e. $-m_1 \vt_1 - m_2
\vt_2 + \arg(u(0)) \in \{ \frac{\pi}{2}, \frac{3\pi}{2} \} $. 
Let $V=V_1 \times V_2 $ be a small open bisector. The
Stokes-directions of $\alpha $ along $D $ in $V $ are defined as
$$
\Sigma_\alpha^{D} := St^{-1} \big( (\frac{\pi}{2}, \frac{3
  \pi}{2}) \big) \subset V \cap \pi^{-1}(D) \ ,
$$
where $St: V \cap \pi^{-1}(D) \to \R/2\pi \Z $ is given by
$$
St(r_1,\vt_1,r_2,\vt_2) := -m_1 \vt_1 - m_2 \vt_2 +
\arg(u \circ \pi(r_1,\vt_1, r_2, \vt_2))
\ .
$$
Remark that $\Sigma^0_\alpha= \Sigma_\alpha^{D} \cap
\pi^{-1}(0) $. We further denote
\begin{equation}\label{eq:Valpha}
V_\alpha := \big( V \smallsetminus \pi^{-1}(D) \big) \cup
\Sigma_\alpha^{D} \ .
\end{equation}
Hence $V_\alpha \cap \pi^{-1}(D) $ consists of those directions in
$\pi^{-1}(D) $ along which $e^{\alpha(x)} $ has rapid
decay for $x $ approaching $D $. Let $\widetilde{\jmath}_{\alpha} : V_\alpha
\hookrightarrow V $ denote the inclusion. Then obviously for $V $
small enough one has
\begin{equation} \label{eq:SDjalph}
\SD|_V = \big( \widetilde{\jmath}_{-\alpha} \big)_! \big( x^{-\lambda}
e^{-\alpha(x)} \cdot \C_{X \smallsetminus D} \big)|_V \ ,
\end{equation}
as well as
$$
\Smod|_V = \big( \widetilde{\jmath}_\alpha \big)_! \big( x^\lambda
e^{\alpha(x)} \cdot \C_{X \smallsetminus D} \big)|_V \ .
$$
Consequently, we deduce an analogous diagram as in \eqref{eq:locV}
where the morphism in the bottom line now reads as
\begin{equation} \label{eq:bottline}
\bR \cHom \big( (\widetilde{\jmath}_{-\alpha})_! \C_{X
  \smallsetminus D} \, , \, \widetilde{\jmath}_! \C_{X \smallsetminus
  D} \big) \lto (\widetilde{\jmath}_{\alpha})_! \, \C_{V_\alpha} \ .
\end{equation}
By the factorization $\widetilde{\jmath} =
\widetilde{\jmath}_{-\alpha} \circ \widetilde{\iota}_{-\alpha} $ with
$\widetilde{\iota}_{-\alpha}: V \smallsetminus \pi^{-1}(D) \hookrightarrow
V_{-\alpha} $, we see that
\begin{multline*}
\bR \cHom \big( (\widetilde{\jmath}_{-\alpha})_!
\C_{V_{-\alpha}} \, , \, \widetilde{\jmath}_! \C_{X \smallsetminus
  D} \big) \cong (\widetilde{\jmath}_{-\alpha})_\ast \bR
\cHom \big( \C_{V_{-\alpha}} ,
(\widetilde{\iota}_{-\alpha})_! \C_{V \smallsetminus \pi^{-1}(D)} \big)  \\
\cong (\widetilde{\jmath}_{-\alpha})_\ast 
\cHom \big( \C_{V_{-\alpha}} ,
(\widetilde{\iota}_{-\alpha})_! \C_{V \smallsetminus \pi^{-1}(D)} \big) =
(\widetilde{\jmath}_{\alpha})_! \, \C_{V_\alpha} \ ,
\end{multline*}
since $(V \smallsetminus V_{\alpha}) \cap \pi^{-1}(D) $ coincides with
the closure of $V_{-\alpha} \cap \pi^{-1}(D) $ inside $\pi^{-1}(D)
$. Hence, \eqref{eq:bottline} is again an isomorphism and thus
$$
\Smod \cong \bR \cHom_{\widetilde{X}} \big( \SD, \widetilde{\jmath}_!
\C_{X \smallsetminus D} \big) 
$$
locally on $\widetilde{X} $ over a crossing point of $D $. Interchanging
$\SD $ and $\Smod $ gives the analogous isomorphism.

In the local situation on $\widetilde{X} $ above a smooth point of $D
$, i.e. above $x=0 $ where locally $D=\{ x_1=0 \} $ on $X $, we proceed in the
same manner, where now $\alpha(x)=x^{-m_1} u(x) $ with $u(0,x_2) \neq
0 $. Instead of bisectors, we consider small tubes $V := V_1 \times
\Delta_2 \subset \widetilde{X} $, where $V_1 $ is a one-dimensional
sector and $\Delta_2 $ a small open disc around $0 $. The definition
of the Stokes-directions of $\alpha $ along $D $ in $V $ are literally
the same as above and the proof of the corresponding isomorphisms
equally holds in this case, completing the proof of the proposition.\\
\qed

\subsection{Rapid decay homology and the asymptotically flat de Rham complex} 

Next, we want to compare the rapid decay complex $\Crdt(\nabla) $
of a flat connection $\nabla $ with the asymptotically flat de
Rham complex in a {\em good} compactification. This will be a key step
in the proof of the main theorem.

\begin{theorem} \label{thm:crddrd}
Let $(X,D) $ be a good compactification of $U $ with respect to the
connection $(E, \nabla) $. We write $d:=\dim_\C(X)=2 $ for its dimension. Then the
rapid decay complex and the asymptotically flat de Rham complex are
isomorphic up to a shift
$$
\Crd_{\widetilde{X}}(\nabla) \cong \DRD(\nabla)[2d] \ .
$$
in the derived category $D^b(\C_{\widetilde{X}}) $.
\end{theorem}
\pf According to Proposition \ref{prop:deg0only}, the right hand side
has cohomology in degree zero only
$$
\SDo \underset{{\rm def}}{=} {\cal
  H}^0(\DRD(\nabla)) \simeq \DRD(\nabla) \ ,
$$
i.e. the inclusion $\SDo \hookrightarrow \DRD(\nabla) $ is a
quasi-isomorphism of complexes.

Now, consider the complex
$\Ctd^{-\bdot} $ of sheaves of
smooth topological chains on $\widetilde{X} $ relative to the
boundary $\widetilde{D} $, defined to be the complex of sheaves
associated to the presheaves
\begin{equation}\label{eq:chainsheaf}
V \mapsto S_{\bdot} \big(
\widetilde{X}, ( \widetilde{X} \smallsetminus V )\cup \widetilde{D}
\big) \ .
\end{equation}
It is standard that this sheaf is a homotopically fine resolution of
the constant sheaf $\C_{\widetilde{X}}[2d] $ shifted by the real
dimension: restricted to $\Uan $ it is the sheaf version of the
absolute singular homology of the real manifold $\Uan $ (cp
\cite{verdier}) and since we take singular homology with respect to
the boundary $\widetilde{D} $, $\Ctd^{-\bdot} $ gives the constant
sheaf along the boundary also: For any small enough $V \subset
\widetilde{X} $, such that $V \subset \Omega $ for a contractible open
$\Omega $ one has by excision
$$
\H^{-p}(\Ctd^{-\bdot}) = H_p \big( \Omega, ( \Omega \smallsetminus V)
\cup \widetilde{D} \big) \cong \left\{ 
  \begin{array}{cl} 0 & \text{for } p \neq 2d \\ \C & \text{for } p=2d
    \ .
\end{array} \right.
$$
We now have a canonical morphism
\begin{equation} \label{eq:ctop}
\Ctd^{-\bdot} \otimes \SDo  \lto \Crdt(\nabla) \ ,
\end{equation}
since for any open $V \subset \widetilde{X} $ and an
asymptotically flat section $\sigma \in \Gamma_V(\SDo) $, by the
characterization of the elements in $\AD $ (see Proposition II.1.1.11
in \cite{sabbah1} or \eqref{eq:ad=rd} above), the section $\sigma \in
\SDo(V) \subset \widetilde{\jmath}_\ast \bve (V) $ is rapidly decaying
along any chain in $c \in \Ctd^{-\bdot}(V) $, hence $c
\otimes \sigma \in \Crdt(\nabla)(V) $. 

We claim that \eqref{eq:ctop} is a quasi-isomorphism. It suffices to
do so for the restriction of the sheaves to a basis of the topology of
$\widetilde{X} $. If $V \subset \widetilde{X} $ is a small open
contractible subset, the image of a $2d $-cell, contained in $\Uan =
\widetilde{X} \smallsetminus \widetilde{D} $, then $\SDo|_V $ is
isomorphic to $\bve|_V $, since $\AD|_{\Uan} \cong \pi^{-1} \O_{\Uan}
$ and on the other hand side 
$$
\H^{-p}( \Crdt(\nabla))(V) \cong H_{p}(\Uan, \Uan \smallsetminus V;
\bve)
$$
coincides with the usual singular homology with values in the local
system $\bve|_{\Uan} $, since no condition of rapid decay is imposed
inside $\Uan $ and for $V \subset \Uan $ we can use an obvious
excision procedure to restrict to the situation inside $\Uan $. This
homology vanishes for $p \neq 2d $ and is isomorphic to $\Gamma(V,
\bve) $ for $p=2d $, the latter isomorphism being induced from the
canonical one 
$$
\H^{-2d}(\Ctd^{-\bdot})(V) \cong H_{2d}(\Uan, \Uan \smallsetminus V)
\cong \C 
$$
for the topological chains only, hence the restriction of
\eqref{eq:ctop} to $\Uan $ is a quasi-isomorphism.

It remains to prove the claim locally around $\widetilde{D} $. Now, since
$(X,D) $ is a good compactification, we can assume that there is a
bicyclic ramification $\rho:Y \to X $ of $D $ with lift $\widetilde{\rho}:\widetilde{Y}
\to \widetilde{X} $ to the oriented real blow-ups such that
locally on $\widetilde{Y} $ the pull-back $\rho^{-1} \nabla $ is
isomorphic to its good formal model. Obviously, the finite map
$\widetilde{\rho} $ induces isomorphisms of complexes
$$
\widetilde{\rho}_\ast \Crd_{\widetilde{Y}}(\rho^{-1} \nabla)
\cong \Crdt(\nabla) \mbox{\quad and \quad} \widetilde{\rho}_\ast {\cal
  C}_{\widetilde{Y}, \widetilde{S}}^{-\bdot}
\cong \Ctd^{-\bdot} \ ,
$$
where $\widetilde{S}= \widetilde{\rho}~^{-1} \widetilde{D} $. Since
$\widetilde{\rho}_\ast \DRdY (\rho^{-1} \nabla) =
\DRD(\nabla) $ (see \eqref{eq:rhoDR} in the proof of Proposition
\ref{prop:deg0only}), and the morphism \eqref{eq:ctop} is obviously
compatible with these isomorphisms, it suffices to prove the assertion on
$\widetilde{Y} $, i.e. we assume $\widetilde{Y}=\widetilde{X}
$, so that for $V \subset \widetilde{X} $ small enough, the connection
is isomorphic on $V $ to its good elementary model and we are reduced
to the case $\nabla=e^\alpha \otimes {\cal R} $ with a regular singular
connection ${\cal R} $.

Now, the solutions of a regular singular connection have moderate
growth and moderate decay at most (see \cite{deligne}), a local
solution of $e^\alpha \otimes {\cal R} $ is asymptotically flat if and
only if $e^\alpha $ has this property. Therefore
$$
\SDo(e^\alpha \otimes {\cal R}) = \SDo(e^\alpha) \otimes
\widetilde{\jmath}_\ast {\cal R}
$$
in the obvious notation for the asymptotically flat solutions
$\SDo(\nabla) $ associated to a given connection.



We proceed as in the proof of Theorem \ref{thm:locdugood}, from
which we also take the notations. First
suppose that $V $ is a small open bisector around some $\vt \in
\pi^{-1}(0) $, where $0 \in D $ denotes a crossing point of $D $. Let
$\Sigma^0_\alpha \subset \pi^{-1}(0) $ denote the set of
Stokes-directions of $e^\alpha $.

If $\vt \in \Sigma^0_\alpha $ we can assume that the directions of the
bisector $V $ are all contained in $\Sigma^0_\alpha $  and then (cp.
\eqref{eq:sdsigma})
$$
\SDo (e^\alpha)|_V = \widetilde{\jmath}_! \big( 
e^{-\alpha} \C_{\Uan} \big) |_V
$$
and consequently, $\SDo(\nabla)|_V = \widetilde{\jmath}_! \bve $. For a
smooth topological chain $c $ in $\widetilde{X} $, a local section $\ve
= e^{-\alpha} \cdot \rho $ in $\bve $ with a local solutions $\rho $
of ${\cal R} $ will not have rapid decay along $c $ in $V $ as
required by the definition unless the chain does not meet
$\widetilde{D} \cap V $. Hence
$$
\Crdt(\nabla)|_V = \Ctd^{-\bdot} \otimes \widetilde{\jmath}_! \bve =
\Ctd^{-\bdot} \otimes \SDo(\nabla)|_V \ .
$$

If $\vt \in \Sigma^0_{-\alpha} $, we can assume that $V $ is an open
bisector such that all the arguments of points in $V $ are contained
in $\Sigma^0_{-\alpha} $. Then
$$
\SDo(e^\alpha)|_V \cong e^{-\alpha} \C_V \ .
$$
Similarly, all twisted chains $c \otimes \ve $ will have rapid
decay inside $V $ and again both complexes considered are equal to
$\Ctd^{-\bdot} \otimes \widetilde{\jmath}_\ast \bve $. 

Finally, if $\vt $ separates the Stokes regions of $\alpha $ and
$-\alpha $, we have \eqref{eq:SDjalph}:
$$
\SDo(e^\alpha)|_V \cong (\widetilde{\jmath}_{-\alpha})_! (e^{-\alpha}
\C_{\Uan})|_V
$$
where $\widetilde{\jmath}_{-\alpha}:V_{-\alpha} \hookrightarrow V $ denotes the
inclusion of the subspace $V_{-\alpha} $ we defined in
\eqref{eq:Valpha}. The characteristic property of $V_{-\alpha} $ is that
$V_{-\alpha} \cap \widetilde{D} $ consists of those directions along
which $e^{-\alpha(x)} $ has rapid decay for $x $ approaching
$\widetilde{D} $. In particular, $c \otimes \ve $ is a rapid decay
chain on $V $ if and only if the topological chain $c $ in $\widetilde{X} $
approaches $\widetilde{D} \cap V $ in $V_{-\alpha} $ at most, i.e.
$$
{\rm im}(c) \cap (\widetilde{D} \cap V) \subset V_{-\alpha} \ .
$$
Again both complexes coincide:
$$
\Crdt(\nabla)|_V = \Ctd^{-\bdot} \otimes (
\widetilde{\jmath}_{-\alpha} )_! \bve =
\Ctd^{-\bdot} \otimes \SDo(\nabla)|_V \ .
$$
The situation for $\vt \in \widetilde{D} $ with $\pi(\vt) $ a smooth
point of $D $ can be handled analogously.

In summary, we have the following composition of quasi-isomorphisms
$$
\SDo[2d] \stackrel{\simeq}{\lto} \C_{\widetilde{X}}[2d] \otimes \SDo
\stackrel{\simeq}{\lto} \Ctd^{-\bdot} \otimes \SDo
\stackrel{\simeq}{\lto} \Crdt(\nabla)
$$
and the claim of the Theorem follows.\\
\qed

\subsection{The general local duality theorem}
\label{subsec:rdundcomp}

Starting with the algebraic connection on
the smooth quasi-projective algebraic surface $U $, we chose a {\em good}
compactification $(X, D) $ (assuming Sabbah's Conjecture holds for the
given connection) in order to define the complex of sheaves
$\Crdt $ on the real oriented blow-up $\widetilde{X} $. By Theorem
\ref{thm:crddrd}, this complex is isomorphic to the
asymptotically flat de Rham complex $\DRD(\nabla^\wee) $ in the
derived category, the proof relying heavily on $X $ being a {\em good}
compactification.

Now, assume that contrary to this situation, we are given an
arbitrary compactification $(X, D) $ with a normal crossing divisor
$D=X \smallsetminus U $ (but without the assumption of a good formal
structure). In this situation, we will still be able to proof the
local duality statement for the de Rham complexes on
$\widetilde{X} $ assuming Sabbah's Conjecture.

To this end, we first study the behavior of the rapid decay and the
moderate de Rham complex under blowing-up a point $z \in D $. More
generally, let $b:Y^{\rm an}
\to \Xan $ be a proper morphism such that $S:=b^{-1}(D) $ is again a normal
crossing divisor and the restriction $b|_{Y^{\rm an} \smallsetminus
  S}: Y^{\rm an} \smallsetminus S \to \Xan \smallsetminus D $ is an isomorphism. Let
$\pi_X:\widetilde{X} \to \Xan $ and $\pi_Y:\widetilde{Y} \to Y^{\rm an} $ denote
the oriented real blow-ups. Lift $b $ to a map
$\widetilde{b}:\widetilde{Y} \to \widetilde{X} $ and consider the open
embeddings $\widetilde{\jmath}:\Xan \smallsetminus D \hookrightarrow
\widetilde{X} $ and $\widetilde{\imath}: Y^{\rm an} \smallsetminus S
\hookrightarrow \widetilde{Y} $.
The de Rham complexes on $\widetilde{X} $ behave well with respect to
this situation, namely
\begin{lemma} \label{lemma:blow}
  There are natural isomorphism
\begin{align*}
\bR\widetilde{b}_\ast (\DRmY(b^\ast \nabla)) & \cong \DRmod(\nabla)
\mbox{ and} \\
\bR\widetilde{b}_\ast (\DRdY(b^\ast \nabla))
& \cong \DRD(\nabla) \ .
\end{align*}
\end{lemma}
\pf Since $\As{X}{? D} $ is flat over $\pi^{-1}(\O_{\Xan}) $, where ? stand
for either $< $ or mod, we conclude by the projection formula that
\begin{multline*}
\bR \widetilde{b}_\ast ( \dr{Y}{? S} ) = \bR \widetilde{b}_\ast \big(
\As{Y}{? S} \Lotimes_{\pi_Y^{-1}(\O_{Y^{\rm an}})} \pi_Y^{-1} {\rm
    DR}_{Y^{\rm an}}(b^\ast \nabla) \big) \cong \\
\cong \bR \widetilde{b}_\ast ( \As{Y}{? S})
\Lotimes_{\pi_X^{-1}(\O_{\Xan})} \pi_X^{-1} \DR(\nabla) \ ,
\end{multline*}
so that it remains to prove $\bR \widetilde{b}_\ast ( \As{Y}{? S} )
\cong \As{X}{? D} $.

Consider the resolution
\begin{equation} \label{eq:resoAd}
\As{Y}{< S} \simeq \big( {\cal P}_{\widetilde{Y}}^{< S}
\otimes_{\pi_Y^{-1}(C^\infty_{Y^{\rm an}})} \pi_Y^{-1} \EYpq{0}{\bdot} ,
\dolb \, \big) \ ,
\end{equation}
where ${\cal P}_{\widetilde{Y}}^{<S} $ denotes the sheaf of $C^\infty
$-functions, which are flat at $\pi_Y^{-1}(S) $. Hence
$$
\bR \widetilde{b}_\ast \As{Y}{<S} 
= \big( \widetilde{b}_\ast {\cal P}_{\widetilde{Y}}^{<S}
\otimes_{\pi_X^{-1}(C^\infty_{\Xan})} \pi_X^{-1} \Epq{0}{\bdot},
\dolb \, \big) \ .
$$
The assertion for $\DRD $ follows, since for any such $b:Y^{\rm an}
\to \Xan $ inducing an isomorphism $Y^{\rm an} \smallsetminus S \to \Xan
\smallsetminus D $, one obviously has
$\widetilde{b}_\ast {\cal P}_{\widetilde{Y}}^{<S} = {\cal
  P}_{\widetilde{X}}^{<D} $.

As for $\DRmod $, the same arguments apply to the resolution
\begin{equation} \label{eq:resoAmod}
\As{Y}{{\rm mod }S} \simeq \big( {\cal P}_{\widetilde{Y}}^{{\rm mod }S}
\otimes_{\pi_Y^{-1}(C^\infty_{Y^{\rm an}})} \pi_Y^{-1} \EYpq{0}{\bdot} ,
\dolb \, \big) \ ,
\end{equation}
with the sheaf ${\cal P}_{\widetilde{Y}}^{{\rm mod }S} $ of $C^\infty
$-functions with moderate growth along $\pi_Y^{-1}(S) $, for which
$\widetilde{b}_\ast {\cal P}_{\widetilde{Y}}^{{\rm mod }S} =
{\cal P}_{\widetilde{X}}^{{\rm mod }D} $ holds as well. \\
\qed

\begin{theorem}[Local Duality] \label{thm:locdual}
  Let $(E, \nabla) $ be an integrable connection on $U $ and $j:U
  \hookrightarrow X $ be an embedding into a smooth projective variety
  such that $D:=X \smallsetminus U $ is a normal crossing
  divisor. Let $\pi:\widetilde{X} \to \Xan $ denote the real oriented
  blow-up of $D $ in $\Xan $. Assuming that Sabbah's Conjecture holds
  for $(E, \nabla) $, the local duality pairing \eqref{eq:locdu} 
$$
\DRmod(E, \nabla) \otimes_\C \DRD(E^\wee, \nabla^\wee) \to \DRD(\O_X,
d)
$$
is non-degenerate.
\end{theorem}
\pf According to Sabbah's Conjecture, there exists a sequence of point blow-ups
$b:Y \to X $ such that $(b^{-1}(E, \nabla))^{\rm an} $ has a good formal
decomposition along $D $ and Proposition \ref{prop:locdugood} can be
applied.

Lifting $b $ to a map between the real oriented blow-ups
$\widetilde{b}: \widetilde{Y} \to \widetilde{X} $, we are in the
situation of Lemma \ref{lemma:blow}. We use the same notation
introduced above and let $W:=Y \smallsetminus S $. Then Proposition
\ref{prop:locdugood} yields an isomorphism
$$
\DRmY(b^\ast \nabla) \stackrel{\cong}{\lto} \bR \cHom_{\widetilde{Y}} \big(
\DRdY(b^\ast \nabla^\wee), \widetilde{\imath}_! \, \C_{W^{\rm an}} \big) \ .
$$
Now $\widetilde{\imath}_! \C_{W^{\rm an}}[2d] $ is the dualizing sheaf
on the manifold $\widetilde{Y} $ with boundary $\widetilde{Y}
\smallsetminus W^{\rm an} $ and similarly $\widetilde{\jmath}_!
\C_{\Uan}[2d] $ for $\widetilde{X} $. 

The resolution \eqref{eq:resoAd} again induces a fine resolution
$$
\widetilde{\imath}_! \C_{W^{\rm an}} \simeq
\big( {\cal P}_{\widetilde{Y}}^{< S}
\otimes_{\pi_Y^{-1}(C^\infty_{Y^{\rm an}})} \pi_Y^{-1}
\EYpq{\bdot}{\bdot} , \partial, \dolb \, \big) \ ,
$$
where the right hand is to be understood as the simple complex
associated to the indicated double complex of Dolbeault-type. Applying
$\widetilde{b}_\ast $ thus yields an isomorphism
\begin{multline*}
\alpha:\bR \widetilde{b}_\ast \widetilde{\imath}_! \C_{W^{\rm an}} =
\widetilde{b}_\ast \big( {\cal P}_{\widetilde{Y}}^{< S}
\otimes_{\pi_Y^{-1}(C^\infty_{Y^{\rm an}})} \pi_Y^{-1}
\EYpq{\bdot}{\bdot} \big) = \\
= \widetilde{b}_\ast {\cal
  P}_{\widetilde{Y}}^{< S} \otimes_{\pi_X^{-1}(C^\infty_{\Xan})}
\pi_X^{-1} \Epq{\bdot}{\bdot} \cong \widetilde{\jmath}_! \C_{\Uan} \ .
\end{multline*}
This construction is the same as the one used for the isomorphisms in
Lemma \ref{lemma:blow} above. It is therefore easy to see that
these isomorphisms are compatible with the local duality pairing given
by the wedge product, i.e. the diagram
\begin{equation} \label{eq:blowwedge}
\begin{CD}
\bR \widetilde{b}_\ast \DRmY (b^\ast \nabla) \otimes_\C \bR
\widetilde{b}_\ast \DRdY(b^\ast \nabla^\wee) @>>> \bR \widetilde{b}_\ast
\widetilde{\imath}_! \C_{W^{\rm an}} \\
@V{\text{Lemma \ref{lemma:blow}}}V{\cong}V @V{\cong}V{\alpha}V \\
\DRmod(\nabla) \otimes_\C \DRD(\nabla^\wee) @>>> \widetilde{\jmath}_!
\C_{\Uan}
\end{CD}
\end{equation}
commutes. Hence, the morphism $\DRmod(\nabla) \to \bR
\cHom_{\widetilde{X}} \big( \DRD(\nabla^\wee), \widetilde{\jmath}_!
\C_{\Uan} \big) $ induced by the lower row of \eqref{eq:blowwedge}
factors as
\begin{equation} \label{eq:rhomwedge}
\begin{CD}
\bR \widetilde{b}_\ast \DRmY (b^\ast \nabla) @>{\beta}>> \bR
\widetilde{b}_\ast \bR \cHom_{\widetilde{Y}} \big( \DRdY(b^\ast
\nabla^\wee), \widetilde{\imath}_! \C_{W^{\rm an}} \big) \\
@V{\cong}VV @VV{\gamma}V \\
\DRmod(\nabla) @>>> \bR \cHom_{\widetilde{X}} \big( \DRD(\nabla^\wee),
\widetilde{\jmath}_! \C_{\Uan} \big) \ ,
\end{CD}
\end{equation}
where $\gamma $ is given by the composition of the natural morphism
$$
\bR \widetilde{b}_\ast \bR \cHom_{\widetilde{Y}} \big( \DRdY,
\widetilde{\imath}_! \C_{W^{\rm an}} \big) \to 
\bR \cHom_{\widetilde{X}} \big( \bR \widetilde{b}_\ast \DRdY, \bR
\widetilde{b}_\ast \widetilde{\imath}_! \C_{W^{\rm an}} \big)
$$
with the morphism $\alpha $ from above. By Poincar\'e-Verdier duality
(Proposition 3.1.10 in \cite{kascha}), $\gamma $ is an isomorphism. In
addition, $\beta $ is an isomorphism due to Theorem
\ref{thm:locdugood}, hence so is the bottom row of
\eqref{eq:rhomwedge}. Interchanging ${\rm DR}^{{\rm mod }D} $ and
${\rm DR}^{<D} $, the same arguments apply, completing the proof of
the theorem.\\ 
\qed

\begin{remark}
The local duality theorem shows that for non-good compactifications
$(X, D) $, the complex $\DRD(\nabla^\wee) $ gives the appropriate dual
object to the moderate (and hence the algebraic) de Rham
cohomology. The resulting pairing, however, lacks
of a similar explicit description we have in case of a good
compactification by period integrals. 

If we assume Sabbah's Conjecture for the given connection, one can
find a finite sequence of point blow-ups
$b:Y \to X $, such that $(Y,S) $ is good with respect to $(E,
\nabla) $. Then $\DRmY(b^\ast \nabla^\wee) $ has cohomology in degree zero
only and is quasi-isomorphic to the rapid decay complex
$\CrdY{\widetilde{Y}} $. On the other hand side, Lemma
\ref{lemma:blow} yields the isomorphism $\DRmod(\nabla) \cong
\bR\widetilde{b}_\ast(\DRmY(b^\ast \nabla)) $, which may have
non-vanishing cohomology in different degrees. In general, we will not
have a topological description of this complex by rapid decay chains
as in the good compactification case.

In the one-dimensional case, every compactification is automatically
good and so the rapid decay complex is always isomorphic to the
asymptotically flat de Rham complex in the derived sense for which the
local duality theorem holds.
\end{remark}

\subsection{Independence of choice of good compactification} \label{sec:indep}

In the definition of the rapid decay homology (Section
\ref{sec:defofrd}) we chose a {\em good} compactification
$(X,D) $ with respect to the given connection $(E, \nabla) $ in order
to define the rapid decay complex on the real oriented blow-up
$\widetilde{X} $. Theorem
\ref{thm:crddrd} and Lemma \ref{lemma:blow} now yield independence of
the choice of the good compactification as an immediate consequence.
\begin{proposition} \label{prop:indep} Assuming Sabbah's Conjecture,
  the following holds: Given two
compactifications $(X_1,D_1) $ and $(X_2,D_2) $ of $U $ such that $D_i=X_i
\smallsetminus U $ are normal crossing
divisors and $(E,\nabla) $ admits a good formal decomposition in each
compactification, there is a canonical isomorphism
$$
\Hy^{-\bdot} ( \widetilde{X}_1, {{\cal C}^{\rm
    rd}_{\widetilde{X}_1}}(\nabla)) \cong \Hy^{-\bdot} (
\widetilde{X}_2, {{\cal C}^{\rm rd}_{\widetilde{X}_2}}(\nabla) ) 
$$
between the corresponding hypercohomologies.
\end{proposition}
\pf We can pass to a common good
compactification $(X,D) $ applying the statement of Sabbah's
Conjecture to the closure of $U $
inside $X_1 \times X_2 $. By Lemma \ref{lemma:blow}, the asymptotically flat
de Rham complexes on $\widetilde{X}_1 $ and $\widetilde{X}_2 $ both
are quasi-isomorphic to the one on $\widetilde{X} $ and hence so are
the rapid decay complexes by Theorem \ref{thm:crddrd}.\\
\qed

\noindent
In particular, the conclusion of the proposition
holds if ${\rm rank} \, E \le 5 $ (by Sabbah's Theorem \ref{thm:sab}).

\section{Global duality of the period pairing}

We will now complete the proof of Theorem \ref{thm:main}, namely the
perfectness of the period pairing for the given flat connection
$(E,\nabla) $ on $U $. 

First, let $(X,D) $ be any compactification of $U $ with a normal
crossing divisor $D $. Taking global sections in the local duality
statement in Theorem \ref{thm:locdual} leads to the following
global duality statement:
\begin{theorem} \label{thm:abstrglob}
Let $(E, \nabla) $ be an integrable connection on the smooth
quasi-projective surface $U $ over $k \subset \C $ and let $j:U \hookrightarrow X $ be an
open embedding into a smooth projective variety $X $ such that $D:=X
\smallsetminus U $ is a normal crossing divisor. Assume that Sabbah's
Conjecture holds for $(E, \nabla) $. Then 
the local duality pairing \eqref{eq:locdu} induces a perfect pairing
\begin{equation}\label{eq:abstrglob}
\big( \Hdr^p(U; E, \nabla) \otimes_k \C \big) \otimes_\C {\mathbb
  H}^{2d-p}(\widetilde{X}, \DRD(\nabla^\wee)) \to \C \ .
\end{equation}
\end{theorem}
\pf If $\widetilde{\jmath}:\Uan \to \widetilde{X} $ denotes the inclusion
of $\Uan $ into the real oriented blow-up of $\Xan $, the local duality
gives an isomorphism
$$
\DRD(\nabla^\wee) \stackrel{\cong}{\lto} \bR \cHom_{\widetilde{X}}(
\DRmod(\nabla), \widetilde{\jmath}_! \C_{\Uan}) \ .
$$
Since $\widetilde{\jmath}_! \C_{\Uan}[2d] $ is the dualizing sheaf on
the compact real $2d $-dimensional manifold $\widetilde{X} $ with boundary
$\widetilde{X} \smallsetminus \Uan $, the local Poincar\'e-Verdier
duality (\cite{kascha}, Prop. 3.1.10) yields the isomorphism
\begin{multline}\label{eq:PV}
\bR \Gamma_{\widetilde{X}} \DRD(\nabla^\wee)[2d] \cong \\
\cong \bR \Gamma_{\widetilde{X}} \bR
\cHom_{\widetilde{X}}( \DRmod, \widetilde{\jmath}_!\C_{\Uan}[2d]) \cong
\Hom_\C^{\bdot} (\bR \Gamma_{\widetilde{X}} \DRmod, \C ) \ ,
\end{multline}
where $\DRmod $ is formed with respect to $\nabla $. Taking $p $-th
cohomology gives the isomorphism we were looking for. Interchanging
$\DRmod $ and $\DRD $, the same arguments prove the perfectness in the
other direction. \\
\qed

To prove our main result, Theorem \ref{thm:main}, we are thus left to
identify the global duality pairing from above with the explicitly given
period pairing in the case of a {\em good} compactification $(X, D)
$. By Proposition \ref{prop:deg0only}, the asymptotically flat de Rham
complex reduces to its $0 $-th cohomology sheaf $\SD $: the inclusion
$$
\SD \hookrightarrow \DRD(\nabla^\wee)
$$
is a quasi-isomorphism. Furthermore, Theorem \ref{thm:crddrd} and the
proof given above establish a quasi-isomorphism of the latter sheaf
shifted by the dimension with the rapid decay complex
$$
\SD[2d] \stackrel{\simeq}{\lto} \SD \otimes \Ctd^{-\bdot} \stackrel{\simeq}{\lto}
\Crdt(\nabla^\wee) \ ,
$$
where $\Ctd^{-\bdot} $ denotes the complex of sheaves of smooth chains
in $\widetilde{X} $ relative to $\widetilde{D} $ (cp. \eqref{eq:chainsheaf}). 

Taking $p $-th cohomology in the duality isomorphism \eqref{eq:PV}
therefore yields the isomorphism 
\begin{equation} \label{eq:PVrd}
\Hy^{-p} (\widetilde{X}, \Crdt(\nabla^\wee)) \stackrel{\cong}{\lto}
\Hom_{\C} \big( \Hy^p(\widetilde{X},\DRmod(\nabla)) , \C
\big) \ .
\end{equation}
We want to prove that \eqref{eq:PVrd} coincides with the morphism
induced by the period pairing.

To this end, consider the resolution
$$
\AD \otimes_{\pi^{-1}\O_{\Xan}} \pi^{-1} \Omega_{\Xan}^r \hookrightarrow
\big( {\cal P}_{\widetilde{X}}^{<D} \otimes_{\pi^{-1} C_{\Xan}^\infty}
\pi^{-1} \Epq{r}{\bdot}, \dolb \big) \ ,
$$
where ${\cal P}_{\widetilde{X}}^{<D} $ as before denotes the sheaf of
$C^\infty $-functions flat at $\pi^{-1}(D) $ and
$\Epq{r}{s} $ denotes the sheaf of $C^\infty $ forms on $\Xan $ of
degree $(r,s) $. This resolution gives rise to the bicomplex
$$
\Rd_{\nabla^\wee}^{\bdot, \bdot}:= \big( {\cal P}_{\widetilde{X}}^{<D}
\otimes_{\pi^{-1}C_{\Xan}^\infty} \pi^{-1} \Epq{\bdot}{\bdot}
\otimes_{\pi^{-1} \O_{\Xan}} \pi^{-1} E^\wee\, , \nabla^\wee, \dolb \big) \
,
$$
the total complex $\Rd_{\nabla^\wee}^{\bdot} $ of which computes
$\SD(\nabla^\wee) \simeq \DRD(\nabla^\wee) $. 
In particular, the complex $\Rd_{(\O_X,d)}^{\bdot} $
associated to the trivial connection $(\O_X, d) $ is a fine resolution of
$\widetilde{\jmath}_! \C_{\Uan} $.

With these resolutions, the local duality pairing $\SD \otimes \Smod
\to \widetilde{\jmath}_! \C_{\Uan} $ (for good compactifications) can
be represented by the bottom row of 
the following obviously commutative diagram:
\begin{equation}\label{eq:finalcommdiag1}
\begin{CD}
\SD[2d] \otimes \Smod @>>> \widetilde{\jmath}_! \C_{\Uan}[2d] \\
@V{\simeq}VV @VV{\simeq}V \\
\Ctd^{-\bdot} \otimes \SD \otimes \Smod @>{\alpha}>> \Ctd^{-\bdot}
\otimes \Rd_{(\O_X, d)}^{\bdot} \ .
\end{CD}
\end{equation}

Recall that $\Dbrd{-\bdot} $ denotes the complex of sheaves of
rapid decay distributions (see section \ref{sec:distri}) on $\widetilde{X} $ which is
a fine resolution of $\widetilde{\jmath}_!\C_{\Uan}[2d] $. There is a
natural quasi-isomorphism
$$
\beta: \Ctd^{-\bdot} \otimes \Rd_{(\O_X,d)}^{\bdot}
\stackrel{\simeq}{\lto} \Dbrd{-\bdot} 
$$
of complexes mapping an element $c \otimes \rho \in \Ctd^{-r}
\otimes \Rd_{(\O_X,d)}^{s}(V)  $ of the left hand
side over some open $V \subset \widetilde{X} $ to the distribution in
$\Dbrd{s-r} $ given by
$\eta \mapsto \int_c \eta \wedge \rho $
for a test form $\eta $ with compact support in $V $. Note hereby,
that the rapid decay property of $\rho $ ensures that the integral
satisfies the estimate \eqref{eq:rddistrib} necessary for the
distribution to be rapidly decaying.

Similarly, we have another natural morphism
$$
\gamma:\Crdt(\nabla^\wee) \otimes \Smod \to \Dbrd{-\bdot} \ , 
\ (c \otimes \ve) \otimes \sigma \mapsto ( \eta \mapsto \int_c \langle
\ve, \sigma \cdot \eta \rangle ) \ .
$$
These morphisms fit into the following diagram
\begin{equation}\label{eq:finalcommdiag2}
\begin{CD}
\Ctd^{-\bdot} \otimes \SD \otimes \Smod @>{\alpha}>> \Ctd^{-\bdot}
\otimes \Rd_{(\O_X, d)}^{\bdot} \\
@V{\simeq}V{\mbox{\scriptsize Thm \ref{thm:crddrd}}}V
@V{\simeq}V{\beta}V \\
\Crdt(\nabla^\wee) \otimes \Smod @>{\gamma}>> \Dbrd{-\bdot} \ .
\end{CD}
\end{equation}
Combining \eqref{eq:finalcommdiag1} and \eqref{eq:finalcommdiag2}, we
see that the period pairing
$$
\Hdr^p(U; E, \nabla) \otimes_\C \Hrd_p(\Uan; E^\wee, \nabla^\wee) \to \C
$$
induced by the bottom row in \eqref{eq:finalcommdiag2} coincides with
the global duality pairing \eqref{eq:abstrglob} induced by the top row
in \eqref{eq:finalcommdiag1} and is therefore perfect by Theorem
\ref{thm:abstrglob}. This completes the proof of our main result,
Theorem \ref{thm:main}.


\bigskip\noindent
\appendix
\section{Existence of solutions with moderate growth}

In this appendix we prove the existence theorem for solutions with
moderate growth for the integrable system of partial differential
equations which arouse in the proof of Proposition
\ref{prop:deg0only}. The proof relies heavily on methods developed by
Hukuhara 
and further elaborated by Ramis, Sibuya (\cite{ramsib}) and
Majima (\cite{majima}). It is a variation of the proof C.~Sabbah gives
for the existence of asymptotically flat solutions (\cite{sabbah2},
Appendix). However, at several places some modifications are
necessary. Note, that we will work entirely in the analytic topology
and therefore omit the superscript 'an' in the following. 

\subsection{The system of differential equations}

Consider the local situation in $X $ around a point $0 \in D $,
where locally $D=\{ x_1x_2=0 \} $ or $D=\{ x_1=0 \} $, and in either
case consider the real oriented blow-up $\pi:\widetilde{X} \to X
$. In order to treat both cases as simultaneously as possible, we use
the notation
$$
\pi: \widetilde{X}= V \times Y \to X \ ,
$$
with $V \cong (\R_+ \times S^1)^n $ and $Y \subset \C^p $ an open disc
around $0 $, where in the case $D=\{ x_1x_2=0 \} $ we have $n=2 $,
$p=0 $ whereas $n=p=1 $ in the case $D=\{ x_1=0 \} $.

\noindent
In the following let $\vt_0 $ be a fixed element in $\pi^{-1}(0)
$ and let there be given
\begin{enumerate}
\item a germ $\alpha \in \Avt{\vt_0} $ on $X$ of the form
  $\alpha(r,\vt):= x_1^{-m_1} x_2^{-m_2} \cdot
  \U(r,\vt) $ with $\U(0,\vt_0) \neq 0 $ in the case $n=2 $ and
  $\alpha(r_1,\vt_1,y) := x_1^{-m_1} \cdot \nu(y) \cdot \U(r_1,\vt_1,y) $
  with $\U(0,\vt_0,y) =1 $ for all $y $ and $\nu(y) $ an invertible
  function.
\item germs $\rho_i \in \Avt{\vt_0} $, $i=1,2 $.
\end{enumerate}
The system of partial differential equations for functions $u $ we
have to consider is the following. For $i=1,2 $, consider the equations
$$
(\Sigma_i): \qquad x_i \frac{\partial u}{\partial x_i} = - x_i
\frac{\partial \alpha}{\partial x_i} \cdot u - \lambda_i u
  + \rho_i =: P_i u + \rho_i \ .
$$
In case $n=2 $, the system $(\Sigma) $ consists of $(\Sigma_1) $ and
$(\Sigma_2) $, whereas in the case $n=p=1 $ we replace $(\Sigma_2) $ by
$$
(\Sigma'): \qquad \frac{\partial u}{\partial y} = \frac{\partial
  \alpha}{\partial y} \cdot u + \rho_2 := P'_2 u +
\rho_2 \ .
$$
The system is called {\em integrable} if $P_1 \rho_2 = P_2 \rho_1 $ in
the first case and $P_1 \rho_2 = P'_2 \rho_1 $ if $n=p=1 $.

\begin{theorem} \label{thm:appendix}
If the system $(\Sigma) $ is integrable and $\rho_i \in \Amodvt{\vt_0}
$ have moderate growth, then there exists a solution $u \in
\Amodvt{\vt_0} $ with moderate growth. A similar statement holds, if
one considers the equation $(\Sigma_1) $ alone.
\end{theorem}

The idea of the proof is to transform the system $(\Sigma) $ into a
fixed point problem for corresponding integral operators. The
technical difficulty, solved by Hukuhara's method, is to choose the
paths of integration in a suitable manner, such that the exponential
part $\exp(\alpha(x)) $ arising from the differential equation is
rapidly decaying whenever the path meets $\pi^{-1}(D) $.

\medskip\noindent
We will now explicitly describe the method in the local case over a
crossing-point of $D $, i.e. $n=2 $ in the above notation. The case of
a smooth point with $n=p=1 $ is the same with some minor modifications
on which we will comment afterwards.

We will have to distinguish two cases, the purely regular and the
irregular case:
\begin{enumerate}
\item The case (IRR): at least one exponent $m_1, m_2 $ is positive,
  say $m_1>0 $.
\item The case (pREG): $m_1=m_2=0 $.
\end{enumerate}

\subsection{The case (IRR)}

\subsubsection{Paths of integration and the fundamental estimate}

We start with introducing some notations. For fixed directions
$\vt=(\vt_1,\vt_2), \eta=(\eta_1, \eta_2) \in \pi^{-1}(0) $ and a
biradius $r=(r_1,r_2) $, which may be $\infty $ also, we write
$$
W(r) := [\vt, \eta]_r := \{ (\rho_1 e^{i \theta_1}, \rho_2 e^{i
  \theta_2}) \in \widetilde{X} \mid 0<\rho_i \le r_i \mbox{ and } \vt_i \le
\theta_i \le \eta_i \} 
$$
for the corresponding closed bisector and $\W $ for its
directions. Given a direction $\vt=(\vt_1, \vt_2) \in \pi^{-1}(0) $,
we say that
\begin{enumerate}
\item $\vt > 0 $, if $\cos(m_1 \vt_1 + m_2 \vt_2 - \arg
  \U(0)) > 0 $; similarly for $<0 $,
\item $\vt $ is of type $(+-) $, if $M(\vt):= m_1 \vt_1 +
  m_2 \vt_2 - \arg \U(0) \equiv \frac{\pi}{2} {\rm mod } 2 \pi $,
\item $\vt $ is of type $(-+) $, if $M(\vt) \equiv \frac{3\pi}{2} {\rm
    mod } 2 \pi $.
\end{enumerate}

For a bisector $W(\infty) $ with directions $\W \subset S^1 \times
S^1 $, we say that $W(\infty)>0 $ or $W(\infty)<0 $ if all $\vt \in \W $
are of this type. We further call $W(\infty) $ is of type $(+-) $ if
$\{ M(\vt) \mid \vt \in \W \} \cap \{ \frac{\pi}{2}, \frac{3 \pi}{2}
\} = \frac{\pi}{2} $, and similar for the type $(-+) $. A bisector is
called {\em proper} if it is of one of these types.

We consider a fixed direction $\vt^0 \in \pi^{-1}(0) $ as in the
theorem. Let $W(\infty)=W_1(\infty) \times W_2(\infty) $ be a proper
bisector around $\vt^0 $, where 
$$
W_1(\infty) = [ \vt_1^0- \delta, \vt_1^0+ \ve] \mbox{ with } \delta, \ve \in
\R_+ \ .
$$
Define
\begin{equation} \label{def:W+-}
W_1(\infty)_+ := \left\{ \begin{array}{cl} [\vt_1^0, \vt_1^0+\ve] &
    \mbox{if } \vt_1^0>0 \\
W_1(\infty) & \mbox{if } \vt_1^0 \mbox{ is of type } (-+)\\
\emptyset & \mbox{otherwise,}
\end{array} \right. 
\end{equation}
$$
W_1(\infty)_- := \left\{ \begin{array}{cl} [\vt_1^0-\delta, \vt_1^0] &
    \mbox{if } \vt_1^0>0 \\
W_1(\infty) & \mbox{if } \vt_1^0 \mbox{ is of type } (+-)\\
\emptyset & \mbox{otherwise.}
\end{array} \right. 
$$
\begin{definition} \label{def:hukfct}
  A {\bf Hukuhara-function} in direction $x_i $ is a piecewise
  affine, decreasing function $a_i: W_i(\infty) \to (0, \pi) $, such that
$$
\cos(m_1 \vt_1+ m_2 \vt_2 - \arg \U(0) + a_i(\vt_i)) \left\{
  \begin{array}{ll} >0 & \mbox{if } \vt_i \in W_i(\infty)_+\\ 
<0 & \mbox{if } \vt_i \in W_i(\infty)_- \ .
\end{array} \right.
$$
\end{definition}
We cite the following lemma from \cite{sabbah2}, Lemme (C.6):
\begin{lemma}
  Any $\vt^0 \in \pi^{-1}(0) $ possesses a basis of neighborhoods
  consisting of proper bisectors admitting Hukuhara-functions in both
  directions.
\end{lemma}

Given such a proper bisector $W(\infty)=[\mu, \eta]_\infty $ with direction
$\W $ together with a Hukuhara-function $a_i $, one defines
$$
A_i: \W_i \to \R_+ \ , \ A_i(\vt_i) := \left\{ \begin{array}{ll} 1 &
    \mbox{if } W_i(\infty)<0 \\
\exp(\Int_{\vt_i}^{\eta_i} \cot(a_i(\tau)) \, d\tau & \mbox{in all
  other cases.}
\end{array} \right.
$$
\begin{definition} \label{def:hukuhara}
The {\bf Hukuhara-domain} inside $W_i(\infty) $ {\bf with
      radius} $r_i $
  is defined as $S_i(r_i) := \{ x_i \mid \arg(x_i) \in \W \mbox{ and }
  0< |x_i| \le r_i \cdot A_i(\arg(x_i)) \} $.
\end{definition}

The domain for our solutions will be
$$
S(r):= \left\{ \begin{array}{ll} S_1(r_1) \times S_2(r_2) & \mbox{if }
    m_2>0 \\[0.2cm]
S_1(r_1) \times W_2(r_2) & \mbox{if } m_2=0 \ , \end{array} \right.
$$
with the standard sector $W_2(r_2) $ with radius $r_2 $.

Remark, that the given fixed $\vt^0 $ endows $S(r) $ with a canonical
base-point $x^0(r)=(x^0_1(r_1), x^0_2(r_2)) $ on the boundary of $S(r)
$ with $\arg(x^0_i(r_i))=\vt^0_i $.

Now, let there be given $r \in \R_+^2 $, $N \in \N^2 $ and $K>0
$. Consider the space $\Fmod(N,r,K) $ of all continuous functions
$\vi:S(r) \to \C $ which are holomorphic in the interior of $S(r) $
and satisfy the estimate
$$
|\vi(x)| \le K \cdot |x|^{-N}:= K \cdot |x_1|^{-N_1} |x_2|^{-N_2} \ .
$$
Consider the vector space $\Fmod(N,r):= \bigcup_{K>0} \Fmod(N,r,K) $,
which we endow with the norm 
$$
|| \vi || := \sup |\vi(x)| \cdot |x|^N \ .
$$

For any $x \in S(r) $, we choose the path $\gamma(t) $ of integration
ending at $x $ depending on the type of the bisector:

\noindent
{\em 1. Case:} $W(\infty)<0 $. Then $\gamma(t)=(\gamma_1(t),x_2) $,
where $\gamma_1 $ is the linear connection between $0 $ and $x_1 $ 
$$
\gamma_1(t):= t x_1/|x_1| \ .
$$

\noindent
{\em 2. Case:} $W(\infty) $ is of type $(+-) $ or $(-+) $. If $\cos(m_1
\arg(x_1)+m_2 \arg(x_2) - \arg(\U(0))) < 0 $, we choose the same path
as in the previous case. Otherwise, choose an auxiliary $\vt_1 \in
\W_1 $ such that $(\vt, \arg(x_2)) $ fulfills this
estimate. The path will then consist of two parts, the first one being
the radial line between $0 $ and $(|x_1|A_1(\vt_1),x_2) $, the second
part, which we will call the {\em Hukuhara-path}, is given by
$$
[\vt_1, \arg{x_1}] \to S_1(r_1) \ , \ \vt \mapsto (|x_1| \exp 
\int_{\vt}^{\arg x_1} \cot (a_1(\tau)) \, dt, 
x_2) \ ,
$$
which we reparameterize by arclength (see the figure below for the
case $(+-)$). 

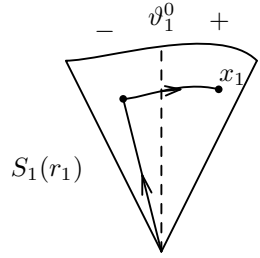
\begin{center}
\begin{texdraw}
\arrowheadtype t:V \arrowheadsize l:0.1 w:0.05
\move (0 0) \lvec (.5 1)
\move (0 0) \lvec (-.5 1)
\clvec (-.3 1)(.3 1.2)(.5 1)
\move (0 0) \avec (-.1 .4) \lvec (-.2 .8) \fcir f:0 r:0.02
\textref h:C v:T \move (0 0)  \lpatt (.05 .05) \lvec (0 1.05)
\move (0 1.3) \htext{$\vt_1^0 $}
\move (-.2 .8) \lpatt() \clvec (-.1 .8)(.1 .9)(.3 .85) 
\fcir f:0 r:0.02
\move (.3 .87) \textref h:L v:B \htext{$x_1$}
\arrowheadsize l:0.08 w:0.05
\move (0 .84) \avec(.1 .85)
\move (-.6 .5) \textref h:C v:T \htext{$S_1(r_1)$}
\move (-.3 1.2) \htext{$-$}
\move (.3 1.25) \htext{$+$}
\end{texdraw}
\end{center}

\noindent
{\em 3. Case:} $W(\infty)>0 $. Assume first that $m_2>0 $. For $i=1,2
$ consider the curve $\gamma_i $ again consisting of two parts, one
radial linear path from the base-point $x^0_i(r_i) $ of $S_i(r_i) $ to
the point $|x_i| A_i(\vt_i^0) $, and another {\em Hukuhara-path} from
this point a path running into $x_i $ parallel to the boundary of
$S_i(r_i) $. The paths over which we will integrate will be of the
form $\gamma(t)=(\gamma_1(t),x_2) $ or $\gamma(t)=(x_1, \gamma_2(t)) $.

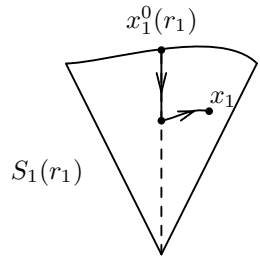
\begin{center}
\begin{texdraw}
\arrowheadtype t:V \arrowheadsize l:0.1 w:0.05
\move (0 0) \lvec (.5 1)
\move (0 0) \lvec (-.5 1)
\clvec (-.3 1)(.3 1.2)(.5 1)
\textref h:C v:T \move (0 0)  \lpatt (.05 .05) \lvec (0 1.05)
\move (0 1.07) \fcir f:0 r:0.02
\move (0 1.3) \htext{$x^0_1(r_1) $}
\move (0 1.05) \lpatt() \avec (0 .85) \lvec (0 .7) \fcir f:0 r:0.02
\clvec (.1 .72)(.15 .77)(.25 .75) 
\fcir f:0 r:0.02
\move (.25 .77) \textref h:L v:B \htext{$x_1$}
\arrowheadsize l:0.08 w:0.05
\move (.15 .745) \avec (.17 .75)
\move (-.6 .5) \textref h:C v:T \htext{$S_1(r_1)$}
\end{texdraw}
\end{center}
In case $m_2=0 $, let $\gamma_1 $ be the same path in $S_1(r_1) $ as
above and let $\gamma_2 $ be the linear path between $x_2^0 $ and $x_2
$.

We emphasize, that all paths are parameterized by arclength. Then one
has the following fundamental estimate (cp. with Lemma 1.7.1 in
\cite{ramsib}):
\begin{lemma} \label{lemma:ineq}
Let $\gamma(t) $ be any of the paths
above of the form $\gamma(t)=(\gamma_1(t),x_2) $ for some $x_2
$. Then, for given $N \in \N^2 $ there exists a pair $r(N) \in
\R_+^2 $ such that for all $r \le r(N) $ and all $x \in S(r) $:
\begin{equation} \label{eq:ineq}
\frac{d}{dt} \big( |\gamma(t)|^{-N} \cdot e^{\Re \alpha(\gamma(t))}
\big) \ge N_1 |\gamma_1(t)|^{-N_1-1} |x_2|^{-N_2} \cdot e^{\Re
  \alpha(\gamma(t))} \ .
\end{equation}
Assuming $m_2>0 $, a similar statement holds for the case
$\gamma(t)=(x_1, \gamma_2(t)) $ interchanging the indices.
\end{lemma}
\pf We have to distinguish the different cases as above.

In case $W(\infty)<0 $, i.e. $e^{\alpha(x)} $ is rapidly decaying in
$W(\infty) $ if $x $ approaches $0 $, $\gamma_1 $ is the radial line
from $0 $ to $x_1 $ and the left hand side of the inequality equals
\begin{equation} \label{eq:radpath}
t^{-N_1-1} |x_2|^{-N_2} \cdot e^{\Re \alpha} \cdot \Big( -N_1+t
\frac{d}{dt} \Re \alpha(\gamma(t)) \Big) \ .
\end{equation}
Now $\alpha\gamma(t)= \gamma_1(t)^{-m_1} x_2^{-m_2}
\U(\alpha\gamma(t)) $ and hence $t \frac{d}{dt}\Re(\alpha\gamma(t))= $
\begin{equation} \label{eq:cosine}
= t^{-m_1} |x_2|^{-m_2} (-m_1 C_1
\cos(m_1 \arg x_1 + m_2 \arg x_2 - \arg \U(\alpha\gamma(t))) + C_2 t)
\end{equation}
where $C_1:= | \U(\alpha\gamma(t))|>0 $ is positive and $C_2:= |
\frac{d}{dt} \U(\alpha\gamma(t)) | \cdot \cos(m_1 \arg x_1 + m_2 \arg x_2- \arg
\U(\alpha\gamma(t))) $ is bounded. Since $W(\infty)<0 $, the cosine 
in \eqref{eq:cosine} is negative for small $r $, i.e. small $t $
(recall that the paths are parameterized by arclength). It follows
that $t \frac{d}{dt} \Re \alpha\gamma(t) $ can be made arbitrarily
large by choosing $r $ small enough.

The same argument applies to the radial part of $\gamma_1 $ in the
case $W(\infty) $ is of type $(+-) $ or $(-+)$, since this radial part
also runs into the rapidly decaying sector inside $S_1(r_1) $. We are
left to consider the Hukuhara-path, which parameterized by the angle $\vt $
reads as $\gamma(\vt):=(\eta(\vt),x_2) $ with
$$
\eta(\vt) = |x_1| \exp  \Big( \int^{\arg(x_1)}_{\vt} \cot (a_1(\tau)) d\tau + i
  \vt \Big)
$$
with the Hukuhara-function $a_1 $ (cp. Definition
\ref{def:hukfct}). 
We thus have 
\begin{equation} \label{eq:dvteta}
\frac{d}{d \vt} |\eta(\vt)| = - |
\eta(\vt)| \cot(a_1(\vt)) \ . 
\end{equation}
Reparametrization by arclength $t $
gives (recall that $0<a_1(\vt)<\pi $)
$$
\frac{dt}{d\vt} = | \frac{d \eta(\vt)}{d \vt} | = | \eta(\vt) (
\cot(a_1( \vt))+i\vt)| = | \eta(\vt)| (\sin(a_1(\vt))^{-1} \ .
$$
It follows that the left hand side of \eqref{eq:ineq}
equals
\begin{equation}\label{eq:sign2}
\pm \Big( N_1 \cos(a_1(\vt))+ \sin(a_1(\vt)) \frac{d}{d\vt} \Re
\alpha(\gamma(\vt)) \Big) \cdot |\gamma(t)|^{-N-1} e^{\Re
  \alpha(\gamma(t))} \ ,
\end{equation}
with the plus sign whenever $\arg(x_1) \in W_1(\infty)_+ $ and the minus sign
otherwise (cp. Definition \ref{def:W+-}), since in the first case the
angle $\vt $ increases along the Hukuhara-path whereas in the last case it
decreases. Now $\Re \alpha(\gamma(\vt))=
|\eta(\vt)|^{-m_1} |x_2|^{-m_2} \cos(\arg(\alpha \gamma(\vt))) $ and
using \eqref{eq:dvteta} one achieves
\begin{multline} \label{eq:sign}
\frac{d}{dt} \Big( |\gamma(t)|^{-N} e^{\Re \alpha\gamma(t)} \Big) = \pm
|\eta(t)|^{-N_1-1} |x_2|^{-N_2} e^{\Re \alpha\gamma(t)} \cdot \\
\cdot \big( -N_1 \cos(a_1(\vt))+m_1 |\eta(\vt)|^{-m_1} |x_2|^{-m_2}
\cdot R \cdot \cos(\rho(\vt)) \big) \ ,
\end{multline}
with $\rho(\vt) $ of the form
$$
\rho(\vt) = m_1 \vt + m_2 \arg(x_2) - \arg
\U(\alpha\gamma(\vt)) +a_1(\vt) \ .
$$
Choosing $r $ small enough, the continuity of $\U $ and the definition
of the Hukuhara-function $a_1(\vt) $ ensure that $\cos(\rho(\vt)) $
has the appropriate sign, namely the same as the sign on the
right hand side of \eqref{eq:sign}. Hence, for $r $ small enough, the
term 
$$
m_1 |\eta(\vt)|^{-m_1} |x_2|^{-m_2} \cdot R \cdot \cos(\rho(\vt))
$$
dominates and can be made arbitrarily large.

The case $W(\infty) >0 $ is similar. The radial
part of the path can be handled in the same way as the radial path in
the case $W(\infty)<0 $, where one has to keep in mind that the
direction is reversed, the path running inwards, and hence the sign in
\eqref{eq:radpath} changes:
$$
\frac{d}{dt} \big( |\gamma(t)|^{-N} \cdot e^{\Re \alpha(\gamma(t))}
\big) = (-1) \cdot \mbox{\eqref{eq:radpath}} \ .
$$
But since $W(\infty)>0 $, the behavior of $\Re \alpha(\gamma(t)) $
changes in the same way, and the same argument as above applies. For the
Hukuhara part, which is contained in either $W_1(\infty)_+ $ or
$W_1(\infty)_- $, the arguments about the sign in \eqref{eq:sign2}
applies and the proof is literally the same as above.\\ 
\qed

\subsubsection{The integral operators} \label{sec:intop}

We now define the integral operators whose fixed points will be the
solutions to the system $(\Sigma) $. We will not go too much into
detail, since -- once the fundamental estimate is established -- this
part of the proof is literally the same as in \cite{sabbah2}, C.13 --
16. Again, the integral operators will differ according to the type of
$W(\infty) $:

If $W(\infty) \not > 0 $, define
$$
Tu(x) := e^{-\alpha(x_1,x_2)} \cdot \Int_{\gamma} (- \lambda_1
u(\xi_1,x_2)+ \rho_1(\xi_1,x_2)) \cdot e^{\alpha(\xi_1,x_2)} \frac{d
  \xi_1}{\xi_1} \
.
$$
For $u \in \Fmod(N,r,K) $, the fundamental estimate yields $|Tu(x)|
\le $
$$
 e^{- \Re \alpha(x)} \cdot \frac{K(1+|\lambda_1|)}{N_1}
|x_2|^{-N_2} \cdot \Int_\gamma \frac{d}{dt} (|\xi_1|^{-N_1} e^{\Re
  \alpha(\xi_1,x_2)} ) dt = \frac{K(1+|\lambda_1|)}{N_1} \cdot |x|^{-N} \ ,
$$
since $e^{\Re \alpha(\xi_1,x_2)} $ decreases rapidly for $\xi_1 \to 0
$. In order to apply the fundamental estimate, we have to choose a $r
\in \R_+^2 $ small enough. Choosing $N $ big enough, such that
$1+|\lambda_1| < N_1 $, we see that $T $ maps $\Fmod(N,r,K) $ into
itself and is contracting. The resulting fixed point $u \in
\Fmod(N,r,K) $ obviously solves $(\Sigma_1) $. Using the integrability
of the system $(\Sigma) $, it follows that $u $ solves the whole
system. We omit the arguments here and refer to \cite{sabbah2}, (B.4)
for the proof. 

If $W(\infty)>0 $, we define $Tu(x):= e^{-\alpha(x)} \cdot (S_1u + S_2
u ) $, where
\begin{equation} \label{eq:defS1}
S_1u:= \int_{\gamma_1} (- \lambda_1 u(\xi_1,x_2) + \rho_1(\xi_1,x_2))
e^{\alpha(\xi_1,x_2)} \frac{d \xi_1}{\xi_1} \mbox{\quad and}
\end{equation}
\begin{equation} \label{eq:defS2}
S_2u:= \int_{\gamma_2} (- \lambda_2 u(x^0_1(r_1),\xi_2) +
\rho_2(x^0_1(r_1),\xi_2)) e^{\alpha(x^0_1(r_1),\xi_2)} \frac{d
  \xi_2}{\xi_2} \ .
\end{equation}
Recall, that $x^0(r)=:x^0 $ denotes the base-point in $S(r) $. Assume
first that $m_2>0 $. Let $N_1,N_2 \ge 2(1+ {\rm max}\{|\lambda_1|,
  |\lambda_2| \}) $. Then the fundamental estimate of Lemma
  \ref{lemma:ineq} yields
$|e^{-\alpha(x^0_1,x_2)} S_2 u | \le $ 
\begin{multline} \label{eq:S2}
\le e^{-\Re
  \alpha(x_0^1,x_2)} \frac{K(1+|\lambda_2|)}{N_2} |x^0_1|^{-N_1}
\int_{\gamma_2} \frac{\partial}{\partial  \xi_2} ( |\xi_2|^{-N_2} e^{\Re
  \alpha(x^0_1,\xi_2)} ) d \xi_2 \le \\
\le \frac{K}{2} |x_1^0|^{-N_1} |x_2|^{-N_2} \ .
\end{multline}
But then
\begin{equation} \label{eq:S2ii}
|e^{-\alpha(x)} S_2u| \le \frac{K}{2} |x_2|^{-N_2} \big( |x_1|^{-N_1}
- e^{\Re \alpha(x)} \int_{\gamma_1} \frac{\partial}{\partial  \xi_1} (
|\xi_1|^{-N_1} e^{\Re \alpha(\xi_1,x_2)} ) d \xi_1 \big) \ .
\end{equation}
Using the fundamental estimate once again gives
\begin{equation} \label{eq:S1}
|e^{-\alpha(x)} S_1 u| \le e^{-\Re \alpha(x)}
\frac{K(1+|\lambda_1|)}{N_1} |x_2|^{-N_2} \int_{\gamma_1}
\frac{\partial}{\partial \xi_1} (
|\xi_1|^{-N_1} e^{\Re \alpha(\xi_1, x_2)} ) d \xi_1 \ ,
\end{equation}
and hence $|Tu| \le \frac{K}{2} |x_1|^{-N_1} |x_2|^{-N_2} $. It
follows that $T $ is a contracting self-map of $\Fmod(N,r,K) $. It is
easy to see that the resulting fixed point $u $ solves the system
$(\Sigma) $, where for the verification of $(\Sigma_2) $ one has to
keep in mind that $(\Sigma) $ is assumed to be integrable.

In the case $W(\infty)>0 $ and $m_2=0 $, the proceeding is
similar. Again, let $Tu=\exp(\alpha(x))(S_1u+S_2u) $ as above. Choose
$r $ small enough, such that
$$
|e^{\alpha(x_1^0,x_2)}| \sup_{\xi_2 \in {\rm im}(\gamma_2)}
|e^{-\alpha(x_1^0, \xi_2)}| < 2 \ ,
$$
which is possible since $m_2=0 $. Then, if $N_2 \ge 4(1+|\lambda_2|) $,
we have
$$
|e^{-a(x_1^0.x_2)} S_2 u| \le 2 K \frac{1+|\lambda_2|}{N_2}
|x_1^0|^{-N_1} ( |x_2|^{-N_2} - |x^0_2|^{-N_2}) \le \frac{K}{2}
|x_1^0|^{-N_1} |x_2|^{-N_2} \ ,
$$
hence \eqref{eq:S2} remains valid. Since $m_1>0 $, the fundamental
estimate on $\gamma_1 $ and then \eqref{eq:S1} are again true, which
proves the claim.

\subsection{The case (pREG) and the case at a smooth point}

The case (pREG) bears no more difficulty, since then
$\alpha(x)=\U(x) $ has no singularity along $\pi^{-1}(D)
$. Replacing the Hukuhara-domain $S(r) $ in the case (IRR) by the
standard bisector $W(r) $ and choosing $\gamma_i $ to be the linear
path between $x_i^0 $ and $x_i $, the integral operator
$$
Tu := e^{-\alpha(x)} (S_1u + S_2 u)
$$
with $S_1 $ and $S_2 $ as in \eqref{eq:defS1} and \eqref{eq:defS2} is
a contracting self-map of $\Fmod(N,r,K) $. The necessary estimates are
easily deduced from the boundedness of $|e^{-\alpha(x)}| $ in that
case.

This completes the proof of Theorem \ref{thm:appendix} in the local
case over a crossing-point of $D $. For a smooth point $0 \in D=\{
x_1=0 \} $, the same arguments apply after the following
replacements in section \ref{sec:intop}. In the case $W(\infty) \not >
0 $, the second integral operator \eqref{eq:defS2} has to be replaced
by
$$
S_2 u := \int_{\gamma_2} \rho_2(x_1^0(r_1), \xi_2)
e^{\alpha(x_1^0(r_1), \xi_2)} \, d \xi_2 \ .
$$
The estimate \eqref{eq:S2} can be achieved again as in the case (pREG)
before, since $|e^{\alpha(x_1^0, x_2)}| $ is bounded for fixed $x_1^0
$. The estimate \eqref{eq:S1} again holds because of the fundamental
estimate in the irregular case $m_1>0 $, or because of the boundedness
of the exponential factor in the purely regular case.

\begin{acknowledgements}
I want to thank Spencer Bloch and H\'el\`ene Esnault for their
permanent interest in this work and for several inspiring discussions.
The present work was started during
  a stay at the University of Chicago supported by the Alexander von
  Humboldt-Stiftung. I thank the Humboldt-Stiftung for financing as well as the
  University of Chicago and especially Spencer Bloch for their
  hospitality. Additionally, I gained from various conversations with
  Stefan Bechtluft-Sachs, Uwe Jannsen, Byungheup Jun and Tomohide
  Terasoma. The present form of this work would not have been possible
  without the comments Claude Sabbah gave to a preliminary version and the
  discussions I had with him during a visit at \'Ecole Polyt\'echnique,
  for which I want to express my deep gratitude.
\end{acknowledgements}

\medskip\noindent
{\footnotesize{\sc NWF I -- Mathematik, Universit\"at Regensburg,
    93040 Regensburg, Germany}\\
marco.hien@mathematik.uni-regensburg.de}


\end{document}